\documentclass[11pt]{article}

\usepackage[margin=1in]{geometry}
\usepackage[authoryear,round]{natbib}

\usepackage[utf8]{inputenc}%
\usepackage[T1]{fontenc}%
\usepackage{hyperref}%
\usepackage{url}%
\usepackage{booktabs}%
\usepackage{amsfonts}%
\usepackage{nicefrac}%
\usepackage{microtype}%
\usepackage{amsmath}
\usepackage{comment}
\usepackage{amssymb}

\usepackage[linesnumbered,ruled]{algorithm2e}
\usepackage{algorithmic}

\usepackage{graphicx}

\newcommand{\papertitle}{Distributed Optimization of \texorpdfstring{$M$}{M}-Convex Functions Subject to Linear Constraints via Continuous Relaxation and Shapley-Folkman Decomposition}

\usepackage{amsthm}
\newtheorem{example}{Example}[section]
\newtheorem{definition}{Definition}[section]
\newtheorem{lemma}{Lemma}[section]
\newtheorem{proposition}{Proposition}[section]
\newtheorem{remark}{Remark}
\newtheorem{theorem}{Theorem}
\newtheorem{corollary}{Corollary}

\title{\papertitle}

\author{%
  Alvaro Maggiar\\
  Amazon.com\\
  \texttt{maggiara@amazon.com}
  \and
  Ozan Candogan\\
  University of Chicago, Booth School of Business\\
  \texttt{ozan.candogan@chicagobooth.edu}
}

\date{}

\begin{document}

\maketitle

\begin{abstract}
$M$-convex functions from discrete convex analysis admit tractable minimization and convex extensions that agree with them at integer points. We study the minimization of a separable collection of $N$ $M$-convex functions subject to $m$ linear coupling constraints. Our two-stage algorithm first solves the convex relaxation and then reconstructs an integer solution through local Carath\'{e}odory decompositions and an LP purification step. The returned solution has no objective increase relative to the convex-relaxation value and a coupling residual that is $O(m)$ for bounded coupling coefficients, with no explicit dependence on $N$. We prove polynomial value-oracle evaluation of the $M$-convex extension and give both an exact value-oracle-polynomial implementation and a distributed consensus-ADMM implementation under a proximal oracle. We also establish proximity and additive envelope-gap guarantees when Stage~1 uses a strongly convex extension. On an inventory placement problem with up to about $10^6$ products, both guarantees hold on every run, and the method is faster than exact integer-programming benchmarks by a factor that grows with the number of regions.

\end{abstract}

\section{Introduction}

$M$-convexity refers to a class of functions from discrete convex analysis (DCA), defined on integer lattice points and motivated as an extension of convexity to the lattice \citep{murota1996convexity,murota2003discrete,murota2022discrete}. Their defining properties, recalled in Section~\ref{sec:mconvexity}, make many $M$-convex integer programs polynomially solvable in standard oracle models despite the NP-hardness of integer programming in general: local optimality implies global optimality, and the lower convex envelope agrees with the function at every integer point.

We study in this paper the problem of minimizing a sum of $N$ $M$-convex functions subject to $m$ linear coupling constraints:
\begin{align}
\min_{x_1,\ldots,x_N}~&\sum_{i=1}^N f_i(x_i)\label{eq:base_problem}\\
\text{s.t.}~& \sum_{i=1}^N A_i x_i = b,\nonumber\\
& x_i\in\mathcal{X}_i\subset\mathbb{Z}^n,~i=1,\ldots,N,\nonumber
\end{align}
where each $f_i:\mathbb{Z}^n\to\mathbb{R}\cup\{+\infty\}$ is $M$-convex with nonempty bounded effective domain $\mathcal{X}_i:=\mathrm{dom}\,f_i$, $A_i\in\mathbb{R}^{m\times n}$, and $b\in\mathbb{R}^m$. For a matrix $A$, $\|A\|_\infty:=\max_k\sum_j|A^{kj}|$ denotes the maximum absolute row sum. We use $\mathbb N:=\mathbb Z_{\ge0}$. Throughout, $x_i\in\mathbb{Z}^n$ is agent $i$'s variable and $x_i^j$ is its $j$-th component. The sets $\mathcal{X}_i$ typically take the form $\mathcal{X}_i=\{x\in\mathbb{N}^n: \sum_{j=1}^n x^j=q_i\}$ for some integer $q_i\in\mathbb{N}$. The assumption that all $\mathcal{X}_i$ have the same ambient dimension $n$ is for ease of presentation, and each could be of dimension $n_i$ instead.

We are especially interested in problems with a large number of functions, e.g. $N$=1M, and low dimensional constraints, e.g. $m$=10. An example of such a problem is given in the placement example (Example~\ref{ex:placement}) below, where we would need to place units from hundreds of thousands of products to a few downstream facilities, subject to capacity constraints.

\paragraph{Contributions.} We make the following contributions.
\begin{enumerate}
    \item \textbf{Bounded-violation reconstruction for sums of $M$-convex functions} (Theorem~\ref{thm:main}). Our two-stage algorithm returns an integer $\bar x\in\prod_i\mathcal{X}_i$ with objective at most the convex-relaxation value and coupling residual at most $m\max_i\|A_i\|_\infty$, with no explicit $N$ factor. The reconstruction uses local decompositions into at most $n$ integer support points and an LP purification step. Appendix~\ref{app:stage1_complexity} gives an exact value-oracle-polynomial implementation, and Appendix~\ref{app:admm_details} an iteration bound for distributed consensus ADMM under a proximal oracle.
    \item \textbf{Polynomial value-oracle evaluation of the $M$-convex convex extension} (Theorem~\ref{thm:poly_extension}). We show that $f^{**}(x^*)$ can be computed with polynomially many value-oracle calls by applying the ellipsoid method to the dual of the Carath\'{e}odory decomposition LP, with an $M$-convex pricing oracle as the separation subroutine.
    \item \textbf{Strongly-convex extension variant.} We analyze strongly convex extensions $\tilde f_i$ that agree with $f_i$ on $\mathcal{X}_i$ and minorize $f_i^{**}$ on $\mathrm{conv}(\mathcal{X}_i)$. We prove proximity and additive envelope-gap bounds while preserving the $\mathcal{O}(m)$ reconstruction residual. Such extensions may then allow for linear-rate Stage~1 solvers.
    \item \textbf{Position in DCA.} The method extends the continuous-relaxation paradigm of \citet{moriguchi2011m} from a single function to linearly coupled separable objectives. It gives a constructive Shapley-Folkman-type guarantee related to \citet{dubois2025frank}, with deterministic LP purification and explicit local hypersimplex supports.
\end{enumerate}

\paragraph{Motivating Example}

\begin{example}[Inventory Placement]\label{ex:placement}
Consider placing $q_i$ units of each of $N$ products across $n$ downstream facilities, where arcs from the source to the facilities may be capacitated. Let $x_i\in\mathbb{Z}^n$ with $x_i^j$ the number of units of product $i$ placed at facility $j$. Typical formulations take $A_i=I$ (unit-based capacities) or $A_i=\mathrm{diag}(w_i)$ (volume-based), and $\mathcal{X}_i=\{x\in\mathbb{N}^n:\sum_j x^j=q_i\}$. The cost $f_i$ represents downstream fulfillment cost per product and is frequently $M$-convex under mild structural assumptions \citep{simchi2014logic,chen2021discrete}.
\end{example}

\subsection{Background}

\subsubsection{\texorpdfstring{$M$}{M}-Convexity}\label{sec:mconvexity}

$M$-convexity was introduced by Murota \citeyearpar{murota1996convexity,murota2003discrete,murota2022discrete} as a lattice analogue of convexity, generalizing the equidistance convexity property. Let $\chi_u\in\{0,1\}^n$ denote the $u$-th unit vector (so $(\chi_u)^v=\delta_{uv}$), write $e_S:=\sum_{u\in S}\chi_u$ for the indicator vector of $S\subseteq\{1,\ldots,n\}$, and for $z,z'\in\mathbb{Z}^n$ let $\mathrm{supp}^+(z-z')=\{j: z^j-z'^j>0\}$ and $\mathrm{supp}^-(z-z')=\{j:z^j-z'^j<0\}$.

\begin{definition}[Exchange Axiom \textbf{M-EXC}($\mathbb{Z}$)]\label{def:exchange}
A function $f:\mathbb{Z}^n\to\mathbb{R}\cup\{+\infty\}$ satisfies \textbf{M-EXC} if for all $x,y\in\mathrm{dom}\,f$ and all $u\in\mathrm{supp}^+(x-y)$ there exists $v\in\mathrm{supp}^-(x-y)$ such that
\begin{align}
    f(x) + f(y) \geq f(x-\chi_u + \chi_v) + f(y+\chi_u - \chi_v).\label{eq:exchange_axiom}
\end{align}
\end{definition}

\begin{definition}[$M$-Convexity]
A function $f:\mathbb{Z}^n\to\mathbb{R}\cup\{+\infty\}$ is \emph{$M$-convex} if it is proper and satisfies \textbf{M-EXC}.
A nonempty set $\mathcal{X}\subset\mathbb{Z}^n$ is \emph{$M$-convex} if its indicator function (0 on $\mathcal{X}$, $+\infty$ outside) is $M$-convex.
\end{definition}

The effective domain $\mathrm{dom}\,f=\{z:f(z)<\infty\}$ of an $M$-convex $f$ lies on a single hyperplane: there exists a constant $q_f$ with $\sum_{j=1}^n z^j = q_f$ for all $z\in\mathrm{dom}\,f$. In particular, any $M$-convex set $\mathcal{X}$ is contained in such a hyperplane; its convex hull $\mathrm{conv}(\mathcal{X})$ is an integral polyhedron whose set of integer points is exactly $\mathcal{X}$, and a polytope as soon as $\mathcal{X}$ is bounded, and a natural canonical form is $\mathcal{X}=\{x\in\mathbb{N}^n:\sum_j x^j=q_f\}$, the integer points of the scaled simplex $q_f\Delta_{n-1}$, which we call the \emph{canonical domain}. Removing the domain-normalization yields the closely related class of $M^\natural$-convex functions; as $M$-convex and $M^\natural$-convex are identical up to dimension lifting, we focus on $M$-convexity.\footnote{This also covers inequality coupling. A rational capacity system $\sum_i A_i x_i\le b$ with nonnegative consumption can first be cleared row by row to obtain integer coefficients and right-hand sides. It can then be written with an integer slack agent. Its box domain is $M^\natural$-convex rather than $M$-convex, so the one-coordinate lift is required to recover fixed-sum form.}

$M$-convexity is tightly connected to substitutability: $M^\natural$-concave set functions on $\{0,1\}^n$ are exactly the gross-substitute valuations \citep{fujishige2003note}, capturing the economic property that raising the price of one item does not decrease demand for others. $M$-convex and $M^\natural$-convex functions appear broadly in economics and game theory \citep{murota2003application,murota2022discrete}, operations \citep{simchi2014logic,chen2021m,chen2021discrete}, and network flows \citep{moriguchi2003capacity,Iwata2005}.

A salient property of $M$-convex functions is local-implies-global optimality, yielding discrete
steepest-descent algorithms \citep{shioura1998minimization,murota2004steepest,oki2025steepest}:
\begin{lemma}[{\citealp[Thm.~6.26]{murota2003discrete}}]\label{lem:local_global}
For $M$-convex $f$ and $x\in\mathrm{dom}\,f$,
\begin{align*}
    f(x)\leq f(y)~\forall y\in\mathbb{Z}^n \Longleftrightarrow f(x) \leq f(x-\chi_u+\chi_v)~\text{for all exchanges}~u,v\in\{1,\ldots,n\}.
\end{align*}
\end{lemma}
The plain iteration count depends on lattice distance and can be pseudopolynomial. Search-space reduction
\citep{shioura1998minimization}, scaling \citep{moriguchi2002scaling}, and long-step variants
\citep{oki2025steepest} provide stronger complexity guarantees.

For a lattice function $f$, let $\bar f:\mathbb{R}^n\to\mathbb{R}\cup\{+\infty\}$ equal $f$ on
$\mathbb{Z}^n$ and $+\infty$ off the lattice. Throughout, $f^{**}$ denotes the real Fenchel biconjugate
of $\bar f$, equivalently its closed lower convex envelope. The key continuous-relaxation bridge is:
\begin{lemma}[{\citealp[Sec.~6.10]{murota2003discrete}}]\label{lem:fstarstar}
For $M$-convex $f:\mathbb{Z}^n\to\mathbb{R}\cup\{+\infty\}$ with effective domain $\mathcal{X}=\mathrm{dom}\,f$, the function $f^{**}$ is convex on $\mathrm{conv}(\mathcal{X})$ and satisfies $f^{**}(x)=f(x)$ for all $x\in\mathcal{X}$.
\end{lemma}
This identity underlies continuous-relaxation algorithms for $M$-convex minimization, most notably the proximity-theorem-based approach of Moriguchi, Shioura, and Tsuchimura \citeyearpar{moriguchi2011m}, which establishes an $L_\infty$ proximity bound of $n-1$ between continuous and discrete minimizers of a single $M$-convex function and uses it to design a two-stage (continuous-then-discrete) algorithm. In~\eqref{eq:base_problem}, the block-separable objective is a direct sum of $M$-convex functions; the difficulty is that a general linear coupling need not preserve the $M$-convex feasible-set structure required by those algorithms.

\paragraph{$M_2$-convex functions and beyond.} In the standard DCA oracle model, the sum of two $M$-convex functions on the same ground set, an $M_2$-convex function, is generally not $M$-convex but remains polynomial-time minimizable via the valuated matroid intersection algorithm of Murota \citeyearpar{murota1996valuatedII,murota1996valuatedI} (see also \citep[Ch.~8]{murota2003discrete}). For three summands, the problem can be shown to be NP-hard.

\subsubsection{Continuous Relaxation of Separable Nonconvex Problems}\label{sec:relaxation_background}

While $M$-convex functions can be continuously extended (Lemma~\ref{lem:fstarstar}), the problem \eqref{eq:base_problem} remains nonconvex due to the integrality of $\mathcal{X}_i$. A standard strategy for large-scale separable nonconvex problems is to solve the convex relaxation and then reconstruct an integer point, either feasible or with a controlled coupling residual.

On the duality-gap side, classical Shapley-Folkman results \citep{starr1969quasi}, together with the
Aubin-Ekeland estimate \citeyearpar{aubin1976estimates} and the Bertsekas-Sandell refinement
\citeyearpar{bertsekas1982estimates}, give $\mathcal{O}(m)$ objective-gap bounds for finite nonconvexity measures
under their compactness and regularity assumptions. Extending a discrete-domain function by $+\infty$
can make that measure infinite. Our result instead permits a controlled coupling residual and compares
the returned objective directly with the convex-relaxation value. Subsequent refinements by Bi-Tang
\citeyearpar{bi2020refined} and Kerdreux-Colin-d'Aspremont \citeyearpar{kerdreux2023stable} tighten classical bounds
under additional assumptions.

On the algorithmic side, \citet{vujanic2014large,vujanic2016decomposition} address MILPs with linear $f_i$,
\citet{udell2014bounding} treats nonconvex $f_i$ with convex $\mathcal{X}_i$, and Dubois-Taine and
d'Aspremont \citeyearpar{dubois2025frank} give a Frank-Wolfe-based two-stage procedure under compactness,
strong-duality, and oracle assumptions. Our construction is closest in spirit to \citet{dubois2025frank}, but specializes to the $M$-convex setting where the Carath\'{e}odory step becomes explicit.

\section{Two-Stage Algorithm and Main Result}\label{sec:algorithm}

\subsection{Overview}\label{sec:overview}

Our approach proceeds in two stages:
\begin{enumerate}
    \item \textbf{Convex relaxation (Stage~1).} Solve the convex relaxation of \eqref{eq:base_problem} using a distributed convex optimization algorithm to obtain continuous primal $x^*=(x_1^*,\ldots,x_N^*)$ with value $f^*_{\mathrm{cvx}}=\sum_i f_i^{**}(x_i^*)$ and dual $\lambda^*$.
    \item \textbf{Integral reconstruction (Stage~2).}
    \begin{enumerate}
        \item First, decompose each $x_i^*$ into a convex combination of integer support points of $\mathcal{X}_i$ via the $M$-convex Carath\'{e}odory lemma (Lemma~\ref{lem:convex_decomposition}).
        \item Adjust the decomposition weights, keeping the coupling satisfied, until all but at most $m$ agents put their whole weight on a single support point, and select one support point per agent from the result (Section~\ref{sec:vertex_selection}). This is the constructive Shapley-Folkman-type step.
    \end{enumerate}
\end{enumerate}
The main result (Theorem~\ref{thm:main}) shows that the returned integer
$\bar x\in\prod_i\mathcal{X}_i$ has objective at most $f^*_{\mathrm{cvx}}$, while permitting a coupling
residual of at most $m\cdot\max_i\|A_i\|_\infty$ and no explicit $N$ factor.

\subsection{Stage 1: Convex Relaxation}\label{sec:stage1}

The convex relaxation of~\eqref{eq:base_problem} replaces each $f_i$ by its biconjugate $f_i^{**}$ and each $\mathcal{X}_i$ by its convex hull:
\begin{align}
\min_{x_1,\ldots,x_N} ~& \sum_{i=1}^N f_i^{**}(x_i) \label{eq:convex_relaxation}\\
\text{s.t.} ~& \sum_{i=1}^N A_i x_i = b, \quad x_i\in\mathrm{conv}(\mathcal{X}_i). \nonumber
\end{align}
$f_i^{**}$ is the lower convex envelope (LCE), and by Lemma~\ref{lem:fstarstar}, $f_i^{**}(x)=f_i(x)$ for all $x\in\mathcal{X}_i$. Problem~\eqref{eq:convex_relaxation} is a convex program, a relaxation of \eqref{eq:base_problem}, and as a result, its optimal value $f^*_{\mathrm{cvx}}$ is a lower bound on \eqref{eq:base_problem}. We assume throughout that~\eqref{eq:convex_relaxation} is feasible; since each $\mathrm{conv}(\mathcal X_i)$ is a polytope and $f_i^{**}$ is finite and continuous on it, it then attains its optimum.

Problem~\eqref{eq:convex_relaxation} has a separable objective and a linear coupling. Our analysis is
agnostic to how it is solved: Stage~1 need only return a point of $\prod_i\mathrm{conv}(\mathcal X_i)$
with small objective excess and coupling residual (Corollary~\ref{cor:eps_stage1}). In particular, the
relaxation can be solved exactly in polynomial time. Under explicit rationality and boundedness
assumptions on the input, Theorem~\ref{thm:exact_stage1} (Appendix~\ref{app:stage1_complexity}) solves
it with polynomially many evaluations of the $f_i$, through an exponentially large configuration LP and
the ellipsoid method. That route settles the complexity question but is centralized and impractical at
scale. For clarity and practicality, we describe here instead a distributed implementation whose
primitive is the proximal oracle
\begin{align}
\mathcal{P}_{i,\rho}(v):=\arg\min_{x_i\in\mathrm{conv}(\mathcal{X}_i)}
\left\{f_i^{**}(x_i)+\frac{\rho}{2}\|x_i-v\|_2^2\right\}.\label{eq:prox_oracle}
\end{align}
Consensus ADMM \citep{boyd2011distributed,maggiar2024consensus} evaluates these $N$ maps in parallel and
projects their aggregate onto the coupling constraint. For exact proximal calls,
Appendix~\ref{app:admm_details} gives the updates and an $O(C_{\mathrm{ADMM}}/\varepsilon)$ ergodic
iteration bound for one-sided objective excess and coupling residual at most $\varepsilon$, hence
$O(NC_{\mathrm{ADMM}}/\varepsilon)$ proximal calls in total, where $C_{\mathrm{ADMM}}$ depends on the
penalty $\rho$ and on the initial distance to an optimal primal-dual pair. The complexity of this route
is therefore stated in proximal calls; Remark~\ref{rem:prox_value_oracle} shows that each call can be
approximated to any prescribed accuracy with polynomially many evaluations of $f_i$. The exact analysis
passes an optimal primal $x^*$ to Stage~2, while Corollary~\ref{cor:eps_stage1} covers a
finite-tolerance output. 

\subsection{Stage 2: Integral Reconstruction}\label{sec:stage2}

The second stage constructs an integer solution $\bar x\in\prod_i\mathcal{X}_i$ from the continuous primal $x^*$; the Stage~1 dual $\lambda^*$ is not needed. It has two components: (i) a per-agent Carath\'{e}odory decomposition of $x_i^*$ into integer support points of $\mathcal{X}_i$; (ii) a vertex-selection step that purifies the decomposition weights, a feasible solution of a small selection LP, until at most $m$ agents are fractional, and then reads one vertex per agent off the result (Section~\ref{sec:vertex_selection}). We analyze each in turn and conclude with the main theorem.

\subsubsection{\texorpdfstring{$M$}{M}-Convex Carath\'{e}odory Decomposition}\label{sec:vertex_decomposition}

\paragraph{Decomposition Existence}

Write $\mathcal{N}(x^*):=\{z\in\mathbb{Z}^n:\lfloor x^{*j}\rfloor\le z^j\le\lceil x^{*j}\rceil\;\forall j\}$ for the integral neighborhood of $x^*$. The following lemma combines the standard local representation of the convex extension of an $M$-convex function with two consequences of the fixed-sum domain: a support bound of $n$ and an explicit hypersimplex form.

\begin{lemma}[$M$-convex Carath\'{e}odory decomposition]\label{lem:convex_decomposition}
Let $f:\mathbb{Z}^n\to\mathbb{R}\cup\{+\infty\}$ be $M$-convex with nonempty bounded effective domain $\mathcal{X}=\mathrm{dom}\,f$, normalized so that $\sum_j x^j=q$ for all $x\in\mathcal{X}$. For any $x^*\in\mathrm{conv}(\mathcal{X})$, there exist integer support points $v_1,\ldots,v_k\in\mathcal{X}$ and weights $\alpha_1,\ldots,\alpha_k\ge 0$ with $\sum_j\alpha_j=1$, $k\le n$, such that:
\begin{align}
x^* &= \sum_{j=1}^k \alpha_j v_j, \label{eq:point_decomposition}\\
f^{**}(x^*) &= \sum_{j=1}^k \alpha_j f(v_j). \label{eq:value_decomposition}
\end{align}
Moreover, writing $I:=\{l:x^{*l}\notin\mathbb{Z}\}$ for the fractional support, each vertex takes the form $v_j=\lfloor x^*\rfloor + e_{S_j}$ with $S_j\subseteq I$ and $|S_j|=p:=\sum_{l=1}^n (x^{*l}-\lfloor x^{*l}\rfloor)$, the \emph{sum of the fractional parts} of $x^*$ and hence the number of coordinates rounded up (an integer, since $\sum_l x^{*l}=q\in\mathbb{Z}$, and $0\le p\le|I|$); in particular, $\|v_j-x^*\|_\infty<1$. The support points thus lie on the translate $\lfloor x^*\rfloor+\{z\in\{0,1\}^n:\sum_lz^l=p\}$ of a \emph{hypersimplex}, the unit cube above $\lfloor x^*\rfloor$ cut by the fixed-sum hyperplane, rather than among the distant vertices $qe_j$ of the canonical domain. There are at most $\binom{|I|}{p}$ feasible candidates, with equality on the canonical domain.
\end{lemma}

\begin{proof}
Every $M$-convex function is $M^\natural$-convex and hence integrally convex
\citep[Thm.~6.42]{murota2003discrete}. Under the bounded-domain assumption, its
local convex extension coincides with its closed lower convex envelope
$f^{**}$ \citep[Sec.~3.4]{murota2003discrete}. Therefore
\begin{align}
f^{**}(x^*)=
\min_{\alpha\ge0}\Big\{
    \sum_{v\in\mathcal X\cap\mathcal N(x^*)}\alpha_v f(v):
    \sum_v\alpha_vv=x^*,\
    \sum_v\alpha_v=1
\Big\}.
\label{eq:decomposition_lp}
\end{align}
This is a feasible finite LP, so it has an optimal basic feasible solution.
Such a solution gives~\eqref{eq:point_decomposition} and
\eqref{eq:value_decomposition}.

Let $B$ denote its equality matrix, whose columns are $(v,1)$. Since
$\sum_{\ell=1}^n v^\ell=q$ for every $v\in\mathcal X$, its rows satisfy
\[
\sum_{\ell=1}^n B_{\ell,\cdot}=qB_{n+1,\cdot}.
\]
Thus $\operatorname{rank}(B)\le n$, and a basic feasible solution has at most
$n$ positive components.

Finally, write $\theta=x^*-\lfloor x^*\rfloor$ and
$I=\{\ell:\theta_\ell>0\}$. Since $\sum_\ell x^{*\ell}=q$,
\[
p:=\sum_\ell\theta_\ell
  =q-\sum_\ell\lfloor x^{*\ell}\rfloor
\]
is an integer. Every $v\in\mathcal X\cap\mathcal N(x^*)$ has the form
$v=\lfloor x^*\rfloor+e_S$ for some $S\subseteq I$. The fixed-sum condition
forces $|S|=p$. Consequently $\|v-x^*\|_\infty<1$, and there are at most
$\binom{|I|}{p}$ candidates. Every such candidate is feasible on the
canonical domain, so equality holds there.
\end{proof}

Applying Lemma~\ref{lem:convex_decomposition} to each agent gives
\begin{align}
x_i^* = \sum_{j=1}^{k_i}\alpha_{ij} v_{ij},\quad f_i^{**}(x_i^*)=\sum_j\alpha_{ij}f_i(v_{ij}),\quad v_{ij}=\lfloor x_i^*\rfloor+e_{S_{ij}},\ |S_{ij}|=p_i.\label{eq:decomposition_i}
\end{align}

Lemma~\ref{lem:convex_decomposition} produces both the integer support points $\{v_{ij}\}$ and the weights
$\{\alpha_{ij}\}$ for each agent. Stage~2(b) receives, for each agent, a \emph{candidate set} containing
these support points and lying in the integral neighborhood of $x_i^*$:
\begin{align}
\{v_{ij}\}_{j=1}^{k_i}\subseteq\mathcal{V}_i
\subseteq\mathcal{X}_i\cap\mathcal{N}(x_i^*).\label{eq:candidate_containment}
\end{align}
Write $\mathcal{V}_i=\{v_{ij}\}_{j=1}^{K_i}$ with $K_i\ge k_i$, and extend
$\alpha_{ij}$ by zero for $j>k_i$.

\paragraph{Computing the decomposition.} The decomposition LP~\eqref{eq:decomposition_lp} has one column
per point of $\mathcal X\cap\mathcal N(x^*)$, and there can be $\binom{|I|}{p}$ of them, exponentially
many in $|I|$. Its pricing problem, however, is $M$-convex and has polynomial value-oracle complexity. We
use this to evaluate the convex extension with polynomially many value-oracle calls, and then to recover
the support points.

To establish polynomial value-oracle complexity of the convex decomposition, we apply the ellipsoid method \citep{khachiyan1980polynomial,bland1981ellipsoid} to the dual of this LP and use the Grötschel--Lovász--Schrijver (GLS) separation-optimization equivalence \citep{grotschel1988geometric}. The GLS theorem states that, for a well-described polyhedron, a polynomial-time strong separation oracle yields polynomial-time strong optimization.

\begin{theorem}[Polynomial value-oracle evaluation of the $M$-convex extension]\label{thm:poly_extension}
Let $f:\mathbb{Z}^n\to\mathbb{R}\cup\{+\infty\}$ be $M$-convex with nonempty bounded effective domain $\mathcal{X}=\mathrm{dom}\,f$, accessible through a value oracle that returns $+\infty$ on $\mathbb{Z}^n\setminus\mathcal{X}$. Assume that $x^*\in\mathrm{conv}(\mathcal{X})\cap\mathbb{Q}^n$, a point $v^0\in\mathcal{X}\cap\mathcal{N}(x^*)$ is supplied, and the finite values of $f$ are rational with a supplied common denominator. After clearing that denominator, assume a supplied integer $\bar W$ bounds $|f(v)|$ for every $v\in\mathcal{X}\cap\mathcal{N}(x^*)$. Then $f^{**}(x^*)$ is computable exactly in
a number of evaluations of $f$ and rational arithmetic operations polynomial in $n$, $\log(1+\bar W)$, and $L_{\mathrm{in}}$, where $L_{\mathrm{in}}$ is the encoding length of $x^*$, $v^0$, and the common denominator.
\end{theorem}

\begin{proof}
We apply the ellipsoid method to the dual of the decomposition LP. The proof proceeds in three steps: formulating the dual, constructing a value-oracle-polynomial separation oracle from $M$-convex minimization, and invoking the Grötschel--Lovász--Schrijver (GLS) theorem.

\emph{Step~1: The dual LP.}
Let $I:=\{j:x^{*j}\notin\mathbb{Z}\}$, $p=\sum_j(x^{*j}-\lfloor x^{*j}\rfloor)$, and $\mathcal{V}:=\{\lfloor x^*\rfloor + e_S : S\subseteq I,\,|S|=p\}\cap\mathcal{X}$, so that $\mathcal{V}=\mathcal{X}\cap\mathcal{N}(x^*)$. The LP dual of the decomposition LP is
\begin{align}
f^{**}(x^*) = \max_{\mu_0\in\mathbb{R},\,\mu\in\mathbb{R}^n}\; \mu_0 + \mu^\top x^* \quad\text{s.t.}\quad \mu_0 + \mu^\top v \le f(v)\;\;\forall\,v\in\mathcal{V}.\label{eq:dec_dual}
\end{align}
For coordinates $j$ with $x^{*j}\in\mathbb{Z}$, every $v\in\mathcal{V}$ has $v^j = x^{*j}$, so $\mu_j$ contributes the same constant $\mu_j x^{*j}$ to both the objective and every constraint; these components can be fixed to zero without loss of generality. The resulting dual has at most $|I|+1\le n+1$ variables. LP strong duality holds since the primal is feasible and bounded.

\emph{Step~2: Separation oracle via $M$-convex minimization.}
Given a candidate $(\mu_0,\mu)$, the separation problem for~\eqref{eq:dec_dual} asks whether $\mu_0 + \mu^\top v \le f(v)$ holds for all $v\in\mathcal{V}$, equivalently whether
\begin{align}
\min_{v\in\mathcal{V}}\big\{f(v) - \mu^\top v\big\} \ge \mu_0.\label{eq:pricing}
\end{align}
If the minimum is less than $\mu_0$, the minimizer $v^*$ provides a violated constraint (a separating hyperplane); otherwise $(\mu_0,\mu)$ is feasible.

The objective $h(v):=f(v)-\mu^\top v$ is $M$-convex, and restricting it to the integer box defining $\mathcal{V}$ preserves $M$-convexity \citep[Thm.~6.13]{murota2003discrete}. Fix a bijection $\phi:I\to\{1,\ldots,|I|\}$. At each iteration choose a steepest feasible exchange and, among ties, lexicographically minimize $(\phi(u),-\phi(v))$. The iteration bound of \citet{minamikawa2021time} then limits the number of iterations to half the $\ell_1$-distance from the starting point to the nearest minimizer, hence to half the $\ell_1$-diameter of the domain, which is at most $p\le n$ here. Each iteration examines $O(n^2)$ exchange pairs, so pricing uses $O(n^3)$ evaluations of $f$. The descent is initialized at the supplied $v^0\in\mathcal{V}$. For the canonical domain, Appendix~\ref{app:offset_sweep} supplies such a point; for a general domain, the phase-one minimization~\eqref{eq:local_phase_one} of Appendix~\ref{app:stage1_complexity} does.

\emph{Step~3: Applying GLS.}
Under the affine shift $\tilde v:=v-\lfloor x^*\rfloor$ and $\tilde\mu_0:=\mu_0+\mu^\top\lfloor x^*\rfloor$, every constraint has $\{0,1\}$ coefficients and an integer right-hand side of magnitude at most $\bar W$. The dual is therefore a well-described rational polyhedron with a strong separation oracle.

For completeness, let $c_v:=(1,\tilde v_I)$ be its constraint normals and let
$\mathcal{L}:=\{d:c_v^\top d=0\ \forall v\in\mathcal{V}\}$ be its lineality space. Writing
$\theta_I:=(x^*-\lfloor x^*\rfloor)_I$, the objective normal $(1,\theta_I)$ belongs to
$\mathrm{conv}\{c_v:v\in\mathcal V\}$, hence to their conic hull. It is therefore orthogonal to
$\mathcal L$, and every recession direction $d$, which satisfies $c_v^\top d\le0$ for all $v$, has
$(1,\theta_I)^\top d\le0$, so none improves the objective. Equivalently, the dual is bounded because the
primal decomposition LP is feasible. Thus any additional affine invariants
of a general $M$-convex domain cause no difficulty. By the strong-separation/strong-optimization
equivalence for well-described polyhedra \citep[Thm.~6.4.9]{grotschel1988geometric}, the exact optimal
value is obtained with a number of strong-separation calls and rational arithmetic operations polynomial
in $n$, $\log(1+\bar W)$, and $L_{\mathrm{in}}$.
Each separation call costs $O(n^3)$ evaluations of $f$, and LP strong duality identifies the returned
value with $f^{**}(x^*)$.
\end{proof}

\begin{remark}[Practical vs.\ theoretical algorithm]\label{rem:practical_vs_theoretical}
The ellipsoid-based algorithm of Theorem~\ref{thm:poly_extension} establishes polynomial value-oracle complexity but is impractical. In the canonical domain, one may solve the decomposition LP by revised simplex with column generation, warm-started by Appendix~\ref{app:offset_sweep}. The restricted master maintains at most $n$ basic columns and is reoptimized after each pricing step. This approach is efficient in practice but has no known polynomial worst-case iteration bound.
\end{remark}

\begin{remark}[The procedure returns both the value and the decomposition]\label{rem:poly_decomposition}
Theorem~\ref{thm:poly_extension} returns the \emph{value} $f^{**}(x^*)$, whereas Stage~2 consumes the
support points. The constructive oracle-dual recovery procedure of
\citet[Lem.~6.5.15]{grotschel1988geometric} makes polynomially many additional separation calls and
returns an optimal dual solution supported on a polynomial-size set of oracle-generated inequalities.
Here each generated inequality, with normal $c_v=(1,\tilde v_I)$, corresponds to the column indexed by
$v$ in the decomposition LP, whose right-hand side is $(1,\theta_I)$ in the shifted coordinates above.
Thus solving the decomposition LP on those columns attains the full optimum.
Choosing a basic restricted optimum and applying the rank argument in
Lemma~\ref{lem:convex_decomposition} leaves at most $n$ positive components. This gives weights $\alpha$
satisfying~\eqref{eq:point_decomposition} and~\eqref{eq:value_decomposition}, supported on points of the
hypersimplex form $\lfloor x^*\rfloor+e_S$, $|S|=p$.
\end{remark}

\paragraph{The Separable Case}

A particular case is that of a separable objective function, which requires no pricing at all, as stated in the proposition below.

\begin{proposition}[Separable costs need no pricing]\label{rem:separable_decomposition}
Let $\mathcal X$ be the canonical domain and $f(x)=\sum_jg_j(x^j)$ with each $g_j$ discrete convex. Write
$\theta:=x^*-\lfloor x^*\rfloor$ and $\sigma_j:=g_j(\lfloor x^{*j}\rfloor+1)-g_j(\lfloor x^{*j}\rfloor)$.
Every $v=\lfloor x^*\rfloor+e_S\in\mathcal X\cap\mathcal N(x^*)$ has
$f(v)=f(\lfloor x^*\rfloor)+\sum_{j\in S}\sigma_j$, which is affine in $e_S$. Hence every decomposition
$x^*=\sum_l\alpha_lv_l$ supported on $\mathcal X\cap\mathcal N(x^*)$ has
\[
\sum_l\alpha_lf(v_l)=f(\lfloor x^*\rfloor)+\sum_j\theta_j\sigma_j=\sum_j\hat g_j(x^{*j}),
\]
where $\hat g_j$ is the chord interpolant of $g_j$, and such a decomposition exists, for instance the
offset sweep of Appendix~\ref{app:offset_sweep}. The function $\sum_j\hat g_j$ is convex, equals $f$ on
$\mathcal X$ and is finite off it, so it is a convex minorant of $\bar f$ and hence at most $f^{**}$ on
$\mathrm{conv}(\mathcal X)$. The display shows it is also at least $f^{**}(x^*)$, since any convex
combination of points of $\mathcal X$ representing $x^*$ bounds $f^{**}(x^*)$ from above. Thus
$f^{**}=\sum_j\hat g_j$, and every local decomposition is value-preserving. In
particular the offset sweep produces the value-preserving
decomposition~\eqref{eq:decomposition_i} in $O(n^2)$ time, with no enumeration, column generation or pricing.
\end{proposition}

\subsubsection{Vertex Selection}\label{sec:vertex_selection}

After decomposition, agent $i$ has a candidate set
$\mathcal V_i=\{v_{ij}\}_{j=1}^{K_i}$ with costs $c_{ij}=f_i(v_{ij})$. Our goal is to select one point
per agent while controlling the aggregate coupling.

\paragraph{LP formulation.}
Consider the LP relaxation of this selection problem, where we let $y_{ij}\ge 0$ be the weight on $v_{ij}$:
\begin{align}
(\mathrm{LP}):~\min_{y} ~& \sum_{i=1}^N \sum_{j=1}^{K_i} c_{ij} y_{ij} \label{eq:lp_selection}\\
\text{s.t.} ~& \sum_{j=1}^{K_i} y_{ij} = 1~\forall i,\quad
\sum_{i=1}^N \sum_{j=1}^{K_i} A_i v_{ij}\, y_{ij} = b,\quad y\ge 0. \nonumber
\end{align}
Setting $y_{ij}=\alpha_{ij}$ from \eqref{eq:decomposition_i} gives a feasible solution with objective $\sum_i f_i^{**}(x_i^*)=f^*_{\mathrm{cvx}}$ (by the value identity \eqref{eq:value_decomposition}); hence, letting $\mathrm{LP}^*$ denote the optimal value of $(\mathrm{LP})$, we have $\mathrm{LP}^*\le f^*_{\mathrm{cvx}}$. The reverse inequality also holds: any feasible $y$ induces $x_i=\sum_j y_{ij}v_{ij}\in\mathrm{conv}(\mathcal{X}_i)$ satisfying the coupling, and convexity of $f_i^{**}$ gives $\sum_j y_{ij}c_{ij}\ge f_i^{**}(x_i)$, so the objective of $(\mathrm{LP})$ is at least the value of the relaxation~\eqref{eq:convex_relaxation}. Hence $\mathrm{LP}^*=f^*_{\mathrm{cvx}}$: over these candidate sets, the selection LP has the same value as the relaxation.

\paragraph{Effective rank.} The number of agents that a solution of $(\mathrm{LP})$ can leave fractional is
governed by
\begin{align}
m_{\mathrm{eff}}:=\dim\mathrm{span}\{A_i(v-v'):v,v'\in\mathcal V_i,\ i=1,\ldots,N\}\le m.\label{eq:m_eff}
\end{align}
If $m=n$ and $A_i=w_iI$, then, since each $\mathcal X_i$ lies in $\sum_jx^j=q_i$, every candidate
difference satisfies $\mathbf 1^\top A_i(v-v')=0$, so $m_{\mathrm{eff}}\le m-1$.

\begin{proposition}[LP sparsity]\label{prop:lp_sparsity}
The constraint matrix of $(\mathrm{LP})$ has rank $N+m_{\mathrm{eff}}$, with $m_{\mathrm{eff}}$ as
in~\eqref{eq:m_eff}. Any basic feasible solution has
at most $N+m_{\mathrm{eff}}$ nonzero $y_{ij}$; hence at most $m_{\mathrm{eff}}$ agents have two or more
nonzero entries, and at least $N-m_{\mathrm{eff}}$ agents have a single $y_{ij^*(i)}=1$.
\end{proposition}

\begin{proof}
A column of $(\mathrm{LP})$ is $(e_i,A_iv_{ij})$. Fix a reference candidate $v_{i1}\in\mathcal V_i$ for each
agent. The $N$ reference columns $(e_i,A_iv_{i1})$ are linearly independent, and every other column
differs from its agent's reference column by $(0,A_i(v_{ij}-v_{i1}))$, a column whose first block is
zero. These differences span the same
space as all pairwise differences, since $A_i(v-v')=A_i(v-v_{i1})-A_i(v'-v_{i1})$. A combination of
reference columns with zero first block has all coefficients zero, so the two spans meet only at $0$ and
the column space has dimension $N+m_{\mathrm{eff}}$. A basic feasible solution has at most as many
nonzero variables as this rank, its support columns being linearly independent. Each agent contributes
at least one nonzero variable, since $\sum_j y_{ij}=1$, so the at most $m_{\mathrm{eff}}$ extra nonzeros
lie in at most $m_{\mathrm{eff}}$ agents, which are the only ones that can be fractional.
\end{proof}

With $m$ in place of $m_{\mathrm{eff}}$, this count is the standard structure of multiple-choice LPs, one
convexity row per block plus $m$ coupling rows, familiar from the multiple-choice knapsack
problem~\citep{sinha1979multiple}.

\paragraph{Purification.} Stage~2 does not need to solve $(\mathrm{LP})$: the decomposition weights
$y=\alpha$ are already feasible with objective $f^*_{\mathrm{cvx}}$, and moving along dependencies among
the fractional agents reduces their number without increasing the cost. For a feasible $y$ write
$J_i(y):=\{j:y_{ij}>0\}$ and $F(y):=\{i:|J_i(y)|\ge2\}$ for its fractional agents.

\begin{proposition}[Purification]\label{prop:purification}
Let $y$ be feasible for~\eqref{eq:lp_selection} with an arbitrary right-hand side $b'$ in place of $b$.
Algorithm~\ref{alg:purification} returns $y'$ that is feasible for the same right-hand side, with
$J_i(y')\subseteq J_i(y)$ for every $i$, $c^\top y'\le c^\top y$, and $|F(y')|\le m_{\mathrm{eff}}$. It
uses at most $\sum_i(|J_i(y)|-1)$ rounds, each computing a null vector of an
$m\times(m_{\mathrm{eff}}+1)$ matrix.
\end{proposition}

\begin{algorithm}[h!]
\SetAlgoNoLine \small
\textbf{Input:} feasible $y$ for~\eqref{eq:lp_selection} with right-hand side $b'$\\[4pt]
\While{$|F(y)|>m_{\mathrm{eff}}$}{
  Choose $m_{\mathrm{eff}}+1$ agents $i\in F(y)$ and, for each, two indices $j_i^+\ne j_i^-$ in $J_i(y)$\;
  Set $u_i\leftarrow A_i(v_{ij_i^+}-v_{ij_i^-})$ and find $\beta\ne0$ with $\sum_i\beta_iu_i=0$\;
  Set $d_{ij_i^+}\leftarrow\beta_i$, $d_{ij_i^-}\leftarrow-\beta_i$, and $d_{ij}\leftarrow0$ otherwise\;
  \lIf{$c^\top d>0$}{$d\leftarrow-d$}
  $t\leftarrow\min\{y_{ij}/(-d_{ij}):d_{ij}<0\}$;\; $y\leftarrow y+t\,d$\;
}
\textbf{Output:} $y$
\caption{Purification of a selection-LP solution}
\label{alg:purification}
\end{algorithm}

\begin{proof}
The vectors $u_i$ lie in a space of dimension $m_{\mathrm{eff}}$, so any $m_{\mathrm{eff}}+1$ of them are
linearly dependent and $\beta$ exists. A chosen agent with $u_i=0$ or $\beta_i=0$ causes no difficulty:
it either receives weight in $\beta$ or contributes no movement. By construction $\sum_jd_{ij}=0$ for every $i$ and
$\sum_i\sum_jA_iv_{ij}d_{ij}=\sum_i\beta_iu_i=0$, so $y+td$ satisfies both equality blocks for every
$t$. Since $\beta\ne0$ and each nonzero $\beta_i$ contributes entries of both signs, $d$ has a negative
entry, all of them on $\bigcup_iJ_i(y)$, so the step $t>0$ is finite and $y+td\ge0$ with at least one
support entry driven to zero. After the sign flip $c^\top d\le0$, so the cost does not increase. No
entry outside the current support becomes positive, and $|J_i|\ge1$ is preserved by
$\sum_jy_{ij}=1$. Each round therefore removes at least one support entry, which bounds the number of
rounds, and the loop exits only when $|F(y)|\le m_{\mathrm{eff}}$.
\end{proof}

When the decomposition supports have $k_i\le n$ points, purifying $\alpha$ takes at most $N(n-1)$ rounds,
each forming $m_{\mathrm{eff}}+1$ vectors $u_i$ in $O(mn)$ operations apiece and computing one null vector
in $O(m^3)$, so the total is linear in $N$. The value $m_{\mathrm{eff}}$ is computed once as a rank, or
replaced by any upper bound $\hat m$ ($m$, or $m-1$ when $m=n$ and $A_i=w_iI$); the algorithm then chooses $\hat m+1$ agents, which are still
dependent, and returns $|F(y')|\le\hat m$. Algorithm~\ref{alg:purification} is the standard reduction
of a feasible LP solution toward a basic one along kernel directions of its support columns (see, e.g.,
the first phase of the proof of \citet[Thm.~2]{megiddo1991finding}, there starting from an optimal
solution), stopped once at most
$m_{\mathrm{eff}}$ agents are fractional. Dubois-Taine and d'Aspremont~\citeyearpar{dubois2025frank} apply the
same reduction, through an exact conic Carath\'{e}odory step after a Frank--Wolfe Stage~1, but carry the
cost as an extra coordinate of the vectors being reduced; this preserves the cost exactly and leaves up to
$m+1$ fractional agents, whereas treating the cost as the LP objective, which only must not increase,
leaves at most $m_{\mathrm{eff}}\le m$. Udell and Boyd~\citeyearpar{udell2014bounding} reach the same
Shapley--Folkman sparsity through a randomized linear objective over the optimal set, and note
purification as the related deterministic route.

\paragraph{Selection rule.} Let $y^\star$ be either the purification of $\alpha$
(Proposition~\ref{prop:purification}) or a basic optimal solution of $(\mathrm{LP})$
(Proposition~\ref{prop:lp_sparsity}). In both cases $y^\star$ is feasible, $c^\top y^\star\le c^\top\alpha$,
and $|F(y^\star)|\le m_{\mathrm{eff}}$. With $J_i:=J_i(y^\star)$, Stage~2 returns
\begin{align}
\bar x_i := v_{i j_i},\qquad j_i\in\arg\min_{j\in J_i} c_{ij},\qquad i=1,\ldots,N.\label{eq:selection_rule}
\end{align}
For at least $N-m_{\mathrm{eff}}$ agents $J_i$ is a singleton and~\eqref{eq:selection_rule} returns its
unique support point.

\paragraph{Constraint violation bound.}

\begin{proposition}[Constraint violation bound]\label{prop:violation_bound}
Let $y$ be feasible for~\eqref{eq:lp_selection} with right-hand side $b'$, write $J_i:=J_i(y)$ and
$F:=F(y)$, and select $\bar x_i=v_{ij_i}$ for any $j_i\in J_i$. With
$D_i:=\max_{v,v'\in\mathcal V_i}\|v-v'\|_\infty$,
\begin{align}
\Big\|\sum_{i=1}^N A_i\bar x_i-b'\Big\|_\infty\le\sum_{i\in F}\|A_i\|_\infty D_i
\le|F|\cdot\max_i\|A_i\|_\infty D_i.\label{eq:violation_bound}
\end{align}
If each $\mathcal V_i\subseteq\mathcal X_i\cap\mathcal N(\xi_i)$ for some $\xi_i$, then $D_i\le1$, and if
moreover $|F|\le m_{\mathrm{eff}}$ the bound is $m_{\mathrm{eff}}\max_i\|A_i\|_\infty\le m\max_i\|A_i\|_\infty$.
In particular, if $m=n$ and $A_i=w_iI$, then $m_{\mathrm{eff}}\le m-1$, as noted after~\eqref{eq:m_eff}, and the bound is at most $(m-1)\max_i|w_i|$, or $m-1$ when $A_i=I$.
\end{proposition}

\begin{proof}
Since $\sum_{j\in J_i}y_{ij}=1$, the rounding error $\delta_i:=\bar x_i-\sum_jy_{ij}v_{ij}$ is the convex
combination $\delta_i=\sum_{j\in J_i}y_{ij}(v_{ij_i}-v_{ij})$, so $\|\delta_i\|_\infty\le D_i$ and
$\|A_i\delta_i\|_\infty\le\|A_i\|_\infty D_i$. For $i\notin F$, $J_i$ is a singleton and $\delta_i=0$.
Since $y$ satisfies the coupling constraint with right-hand side $b'$, the residual equals
$\sum_{i\in F}A_i\delta_i$, which proves~\eqref{eq:violation_bound}.
\end{proof}

\subsubsection{Main Theorem}\label{sec:main_theorem}

We now combine the Stage-2 components into a single performance guarantee.

\begin{theorem}[Main result]\label{thm:main}
Suppose Stage~1 returns an optimal $x^*$ for \eqref{eq:convex_relaxation} with value
$f^*_{\mathrm{cvx}}$. For each agent, let the candidate set satisfy
$\mathcal{V}_i\subseteq\mathcal{X}_i\cap\mathcal{N}(x_i^*)$ and contain the support of a value-preserving decomposition~\eqref{eq:decomposition_i}. Let $y^\star$ be chosen as in the selection rule, for instance the purification of $\alpha$, and let $\bar x$ be the selection~\eqref{eq:selection_rule}. Then $\bar x\in\prod_i\mathcal{X}_i$ and
\begin{align}
\sum_{i=1}^N f_i(\bar{x}_i) &\le f^*_{\mathrm{cvx}},\label{eq:main_obj_bound}\\
\Big\|\sum_{i=1}^N A_i\bar{x}_i - b\Big\|_\infty &\le m\cdot\max_i\|A_i\|_\infty.\label{eq:main_viol_bound}
\end{align}
\end{theorem}

\begin{proof}
By Lemma~\ref{lem:convex_decomposition}, $\sum_i f_i^{**}(x_i^*)=\sum_i\sum_j\alpha_{ij}f_i(v_{ij})=f^*_{\mathrm{cvx}}$, so $y=\alpha$ is feasible for the selection LP~\eqref{eq:lp_selection} with $c^\top\alpha=f^*_{\mathrm{cvx}}$. For the selection~\eqref{eq:selection_rule}, every agent satisfies $f_i(\bar x_i)=\min_{j\in J_i}c_{ij}\le\sum_j y^\star_{ij}c_{ij}$, since a minimum over $J_i$ is at most any convex combination of $\{c_{ij}\}_{j\in J_i}$ (with equality when $J_i$ is a singleton). Summing over $i$,
\begin{align*}
\sum_i f_i(\bar x_i)\;\le\;c^\top y^\star\;\le\;c^\top\alpha\;=\;f^*_{\mathrm{cvx}},
\end{align*}
which is~\eqref{eq:main_obj_bound}. The constraint residual~\eqref{eq:main_viol_bound} is
Proposition~\ref{prop:violation_bound} with $\xi_i=x_i^*$, which in fact gives the sharper
$m_{\mathrm{eff}}\cdot\max_i\|A_i\|_\infty$.
\end{proof}

\begin{corollary}[Inexact or coupling-infeasible Stage~1]\label{cor:inexact_stage1}\label{cor:eps_stage1}
Let $x^*\in\prod_i\mathrm{conv}(\mathcal{X}_i)$ be any Stage~1 point, neither assumed optimal
for~\eqref{eq:convex_relaxation} nor assumed to satisfy its coupling constraint, and set
\begin{align}
\hat b:=\sum_{i=1}^N A_i x_i^*,\qquad r:=\hat b-b.\label{eq:shifted_rhs}
\end{align}
For each agent, let $\mathcal{V}_i\subseteq\mathcal{X}_i\cap\mathcal{N}(x_i^*)$ contain the support of a
value-preserving decomposition of $x_i^*$. Let $y^\star$ be the purification of $\alpha$ for the
selection LP~\eqref{eq:lp_selection} with right-hand side $\hat b$, or a basic optimal solution of that LP,
and let $\bar x$ be the selection~\eqref{eq:selection_rule}. Then
$\bar x\in\prod_i\mathcal{X}_i$ and
\begin{align}
\sum_{i=1}^N f_i(\bar x_i)&\le\sum_{i=1}^N f_i^{**}(x_i^*),\label{eq:inexact_obj_bound}\\
\Big\|\sum_{i=1}^N A_i\bar x_i-b\Big\|_\infty&\le m_{\mathrm{eff}}\cdot\max_i\|A_i\|_\infty+\|r\|_\infty,
\label{eq:inexact_viol_bound}
\end{align}
with $m_{\mathrm{eff}}\le m$ as in~\eqref{eq:m_eff}, computed for these candidate sets.
If, for some $\varepsilon_{\mathrm{obj}},\varepsilon_{\mathrm{feas}}\ge0$,
\begin{align}
\sum_{i=1}^N f_i^{**}(x_i^*)&\le f^*_{\mathrm{cvx}}+\varepsilon_{\mathrm{obj}},\nonumber\\
\|r\|_\infty&\le\varepsilon_{\mathrm{feas}},\label{eq:stage1_tolerances}
\end{align}
then
\begin{align}
\sum_{i=1}^N f_i(\bar x_i)&\le f^*_{\mathrm{cvx}}+\varepsilon_{\mathrm{obj}},\label{eq:eps_obj_bound}\\
\Big\|\sum_{i=1}^N A_i\bar x_i-b\Big\|_\infty&\le
m_{\mathrm{eff}}\cdot\max_i\|A_i\|_\infty+\varepsilon_{\mathrm{feas}}.
\label{eq:eps_viol_bound}
\end{align}
If $x^*$ is optimal for~\eqref{eq:convex_relaxation} then $r=0$ and $\sum_i f_i^{**}(x_i^*)=f^*_{\mathrm{cvx}}$,
recovering Theorem~\ref{thm:main}.
\end{corollary}

\begin{proof}
Lemma~\ref{lem:convex_decomposition} uses only $x_i^*\in\mathrm{conv}(\mathcal{X}_i)$ and not the coupling
constraint, so it applies verbatim to each agent and yields~\eqref{eq:decomposition_i}. Setting
$y_{ij}=\alpha_{ij}$ gives $\sum_i\sum_j A_iv_{ij}y_{ij}=\sum_i A_ix_i^*=\hat b$ by~\eqref{eq:shifted_rhs}, so
$y=\alpha$ is feasible for the shifted LP with cost
$c^\top\alpha=\sum_i\sum_j\alpha_{ij}f_i(v_{ij})=\sum_i f_i^{**}(x_i^*)$ by the value
identity~\eqref{eq:value_decomposition}. Propositions~\ref{prop:purification} and~\ref{prop:lp_sparsity}
hold for any right-hand side, so $y^\star$ is feasible for $\hat b$, has $c^\top y^\star\le c^\top\alpha$,
and has at most $m_{\mathrm{eff}}$ agents with $|J_i|>1$. The argument of Theorem~\ref{thm:main} then gives
$\sum_i f_i(\bar x_i)\le c^\top y^\star\le\sum_i f_i^{**}(x_i^*)$,
which is~\eqref{eq:inexact_obj_bound}. Proposition~\ref{prop:violation_bound} applied with right-hand side
$b'=\hat b$ and $\xi_i=x_i^*$ gives $\|\sum_i A_i\bar x_i-\hat b\|_\infty\le m_{\mathrm{eff}}\max_i\|A_i\|_\infty$,
and~\eqref{eq:inexact_viol_bound} follows by the triangle inequality with $\hat b-b=r$.
Under~\eqref{eq:stage1_tolerances}, substituting the two tolerance bounds into
\eqref{eq:inexact_obj_bound} and~\eqref{eq:inexact_viol_bound} gives
\eqref{eq:eps_obj_bound} and~\eqref{eq:eps_viol_bound}.
\end{proof}

\begin{remark}[Inexact or surrogate Stage~1]\label{rem:inexact_stage1}
If $x^*$ violates the coupling constraint, the selection LP with right-hand side $b$ need not be feasible.
The corollary therefore recenters it at $\hat b=\sum_iA_ix_i^*$. For consensus ADMM, Stage~2 uses the
$x$-iterate or its ergodic average, which lies in $\prod_i\mathrm{conv}(\mathcal X_i)$, and charges its
coupling residual through $\varepsilon_{\mathrm{feas}}$.
\end{remark}

\begin{remark}[What the tolerances require of the solver]\label{rem:solver_tolerances}
For consensus ADMM,~\eqref{eq:admm_ergodic_bounds} gives~\eqref{eq:stage1_tolerances} after
$K\ge C_{\mathrm{ADMM}}\max\{\varepsilon_{\mathrm{obj}}^{-1},
\varepsilon_{\mathrm{feas}}^{-1}\}$ iterations.
The coupling residual is directly computable. For the objective, weak duality gives
$\lambda^\top b+\sum_i\min_{x_i\in\mathcal X_i}\{f_i(x_i)-\lambda^\top A_ix_i\}\le f^*_{\mathrm{cvx}}$
for every $\lambda\in\mathbb R^m$, where each term is an $M$-convex minimization, so
$\sum_if_i^{**}(x_i^*)$ minus this quantity, evaluated at the Stage~1 multiplier, is a computable
upper bound on the objective excess. The feasibility
tolerance adds to a residual bound that is already $O(m_{\mathrm{eff}}\max_i\|A_i\|_\infty)$, so for
$\varepsilon_{\mathrm{feas}}\le\max_i\|A_i\|_\infty$ the Shapley-Folkman term dominates
and~\eqref{eq:eps_viol_bound} is within a factor $1+1/m$ of~\eqref{eq:main_viol_bound}.
Appendix~\ref{app:stage1_complexity} gives an independent exact value-oracle-polynomial implementation
of Stage~1, and sketches an approximate value-oracle implementation of each proximal call.
\end{remark}

\begin{remark}\label{rem:integral_convexity}
A natural question is to ask where $M$-convexity is required, as opposed to mere integral convexity \citep{favati1990convexity}, and whether the same scheme could be applied to more general integral convex functions. The formal convex relaxation and reconstruction argument extend, but the polynomial value-oracle implementations do not. The $M$-convex structure is used in two ways:
\begin{itemize}
    \item \textbf{The Carath\'{e}odory decomposition} (Lemma~\ref{lem:convex_decomposition}) follows from the lower-convex-envelope representation and integral convexity, with at most $n+1$ integer support points in the integral-convexity simplex \citep{favati1990convexity} containing $x^*$ (instead of $n$ for the $M$-convex case due to the domain normalization constraint). However, the explicit hypersimplex form $v_j=\lfloor x^*\rfloor+e_{S_j}$, $|S_j|=p$, is specific to $M$-convex sets via their domain normalization $\sum_j x^j=q$, and does not carry over.
    \item \textbf{The polynomial value-oracle evaluation of the convex extension} (Theorem~\ref{thm:poly_extension}) relies on $M$-convexity through the pricing subproblem. $\min_{v\in\mathcal{V}}\{f(v)-\mu^\top v\}$ is $M$-convex because adding a linear term preserves the exchange axiom, and its $O(n^3)$ value-oracle bound follows from $M$-convex steepest descent on the unit-diameter hypersimplex slice. No polynomial value-oracle analogue is possible for arbitrary integrally convex functions without further structure.%
\end{itemize}
\end{remark}

\subsubsection{Candidate Sets Without Weights and the Selection Dual}\label{sec:selection_dual}

\paragraph{Full enumeration.}
Write $I_i:=\{j:x_i^{*j}\notin\mathbb{Z}\}$ for agent $i$'s fractional support. When $\binom{|I_i|}{p_i}$ is tractable ($|I_i|$ small, or $p_i$ close to $0$ or $|I_i|$), one may take $\mathcal V_i$ to be \emph{all} vertices $\{\lfloor x_i^*\rfloor+e_S:S\subseteq I_i,\,|S|=p_i\}$ and skip the decomposition computation entirely. For the canonical domain of Section~\ref{sec:mconvexity} every such vertex lies in $\mathcal{X}_i$; for a general $M$-convex domain the set must first be intersected with $\mathcal{X}_i$ through the value oracle, and the inclusion~\eqref{eq:candidate_containment} survives that filtering because the support points $v_{ij}$ of Lemma~\ref{lem:convex_decomposition} are themselves feasible. Enumeration returns the points but not the weights $\alpha$.

\paragraph{Recovering a selection.} Without weights, $y=\alpha$ is not available to purify. Instead, solve
the selection LP to optimality, with right-hand side $b$ for Theorem~\ref{thm:main} or
$\hat b=\sum_iA_ix_i^*$ for Corollary~\ref{cor:eps_stage1}, and take a basic optimal $y^\star$. The weights
$\alpha$ of Lemma~\ref{lem:convex_decomposition} exist even though they were not computed, and they are
feasible for that LP, so $c^\top y^\star=\mathrm{LP}^*\le c^\top\alpha$. Together with
Proposition~\ref{prop:lp_sparsity}, $y^\star$ satisfies the selection rule, and both results hold unchanged.

\paragraph{The selection dual.} At scale, $(\mathrm{LP})$ can be approached through its Lagrangian dual.
Dualizing only the coupling constraint, with the sign convention of Section~\ref{sec:stage1}, gives
\begin{align}
L(\lambda) = \lambda^T b + \sum_{i=1}^N \min_{j\in\{1,\ldots,K_i\}}
\big\{ c_{ij} - \lambda^T A_i v_{ij} \big\}. \label{eq:lagrangian}
\end{align}
The per-agent subproblem minimizes a linear function over the simplex $\{y_i\ge 0:\sum_j y_{ij}=1\}$ and its minimum is always attained at a vertex $y_{ij^*(i,\lambda)}=1$. Hence every Lagrangian iterate is an integer selection, and $\max_\lambda L(\lambda)=\mathrm{LP}^*$ by LP strong duality. The dual can be maximized by subgradient ascent,
\begin{align}
\lambda_{t+1} = \lambda_t + s_t\Big(b - \sum_{i=1}^N A_i \bar{x}_{i,t}\Big),\quad \bar{x}_{i,t}=v_{i,j^*(i,\lambda_t)},\label{eq:subgradient}
\end{align}
with $s_t$ a diminishing step, or a Polyak step with target $f^*_{\mathrm{cvx}}$, which is exact because
$\mathrm{LP}^*=f^*_{\mathrm{cvx}}$ whenever the candidate sets contain a value-preserving decomposition,
as full enumeration guarantees. A subgradient iterate is an integer selection but carries no residual
guarantee: a feasible primal $y$ must still be recovered before Proposition~\ref{prop:violation_bound}
applies.

\begin{remark}[Reuse of the Stage~1 multiplier]\label{rem:multiplier_reuse}
If Stage~1 also returns an optimal multiplier
$\lambda^*$, then each $x_i^*$ minimizes
$f_i^{**}(x_i)-(\lambda^*)^\top A_ix_i$. Every candidate $v_{ij}\in\mathcal X_i$ therefore has
$c_{ij}-(\lambda^*)^\top A_iv_{ij}$ at least that minimum, with equality on the support of the
value-preserving decomposition. Hence $L(\lambda^*)=f^*_{\mathrm{cvx}}=\mathrm{LP}^*$: the multiplier
is an exact dual certificate. It is not needed for the selection itself, which purifies $\alpha$, and,
being already dual-optimal, it is not a useful starting point for~\eqref{eq:subgradient}: a Polyak step
from it is zero.
\end{remark}

\section{Strongly Convex Extensions}\label{sec:strongly_convex_extension}

\subsection{Motivation}

A strongly convex extension can make the Stage~1 subproblems simpler than the piecewise-linear
lower convex envelope, and may allow faster Stage~1 solvers. Let $\tilde f_i$ agree with $f_i$ on
$\mathcal X_i$ and minorize $f_i^{**}$ on $\mathrm{conv}(\mathcal X_i)$. Stage~1 then returns the minimizer
of the surrogate rather than of the envelope, so Stage~2 decomposes a different point. This section bounds
the distance between the two minimizers and the resulting loss in value, and identifies when the exact
bound of Theorem~\ref{thm:main} survives.

\subsection{Proximity Bounds for Strongly Convex Extensions}\label{sec:proximity_bounds}

Let $G(x):=\sum_i f_i^{**}(x_i)$ and $\tilde G(x):=\sum_i\tilde f_i(x_i)$ denote the LCE and surrogate aggregate objectives. Throughout this section, $\tilde G$ is $\mu$-strongly convex on $\prod_i\mathrm{conv}(\mathcal{X}_i)$ and, for every agent,
\begin{align*}
\tilde f_i(v)=f_i(v)\quad(v\in\mathcal X_i),\qquad
\tilde f_i(x_i)\le f_i^{**}(x_i)\quad(x_i\in\mathrm{conv}(\mathcal X_i)).
\end{align*}
Define the per-agent and total envelope gaps
\begin{align}
\Delta_i:=\sup_{x_i\in\mathrm{conv}(\mathcal{X}_i)}\big[f_i^{**}(x_i)-\tilde f_i(x_i)\big],\quad \Delta:=\sum_{i=1}^N\Delta_i.\label{eq:envelope_gap_total}
\end{align}

\begin{proposition}[Proximity bound]\label{prop:proximity}
Let $\mathcal{F}:=\big\{x\in\textstyle\prod_i\mathrm{conv}(\mathcal{X}_i):\sum_i A_i x_i=b\big\}$ denote the convex relaxation feasible set. Let $x^*\in\arg\min_{\mathcal{F}}G$ be a minimizer of the LCE objective, and let $\tilde x^*\in\arg\min_{\mathcal{F}}\tilde G$ be the (unique, by strong convexity) minimizer of the surrogate. Then
\begin{align}
\|\tilde x^*-x^*\|_2\;\le\;\sqrt{\frac{2\Delta}{\mu}}.\label{eq:proximity}
\end{align}
\end{proposition}
\begin{proof}
See Appendix~\ref{app:proximity_proofs}.
\end{proof}

\begin{corollary}[Additive value guarantee at the surrogate minimizer]\label{cor:additive_value}
Under the hypotheses of Proposition~\ref{prop:proximity}, run Stage~2 at $\tilde x^*$ and suppose each candidate set satisfies $\mathcal{V}_i\subseteq\mathcal{X}_i\cap\mathcal{N}(\tilde x_i^*)$ and contains the vertices of an LCE decomposition of $\tilde x_i^*$, that is, of a family satisfying~\eqref{eq:point_decomposition}--\eqref{eq:value_decomposition} at $\tilde x_i^*$. Then the selection~\eqref{eq:selection_rule} satisfies $\bar x\in\prod_i\mathcal{X}_i$,
\begin{align}
\sum_{i=1}^N f_i(\bar x_i)\;\le\;G(\tilde x^*)\;\le\;f^*_{\mathrm{cvx}}+\Delta,\label{eq:additive_value_gap}
\end{align}
with the residual bound~\eqref{eq:main_viol_bound} of Theorem~\ref{thm:main} unchanged.
\end{corollary}
\begin{proof}
The assumed inclusion makes $y_{ij}=\alpha_{ij}$ feasible for the selection LP~\eqref{eq:lp_selection} at $\tilde x^*$, which is coupling-feasible since $\tilde x^*\in\mathcal{F}$, with objective $\sum_i f_i^{**}(\tilde x_i^*)=G(\tilde x^*)$. The selection step of the proof of Theorem~\ref{thm:main} gives $\sum_i f_i(\bar x_i)\le c^\top y^\star\le c^\top\alpha=G(\tilde x^*)$. Moreover,
$G(\tilde x^*)\le\tilde G(\tilde x^*)+\Delta\le\tilde G(x^*)+\Delta\le G(x^*)+\Delta
=f^*_{\mathrm{cvx}}+\Delta$. The residual bound follows from
Proposition~\ref{prop:violation_bound} with $\xi_i=\tilde x_i^*$.
\end{proof}

\begin{corollary}[Inexact surrogate Stage~1]\label{cor:inexact_surrogate}
Let $\tilde f^*_{\mathrm{cvx}}:=\min_{x\in\mathcal F}\tilde G(x)$, and suppose
$\hat x\in\prod_i\mathrm{conv}(\mathcal X_i)$ satisfies
\begin{align*}
\tilde G(\hat x)\le\tilde f^*_{\mathrm{cvx}}+\varepsilon_{\mathrm{obj}},\qquad
\Big\|\sum_iA_i\hat x_i-b\Big\|_\infty\le\varepsilon_{\mathrm{feas}}.
\end{align*}
If the candidate sets lie in $\mathcal X_i\cap\mathcal N(\hat x_i)$ and contain value-preserving LCE
decompositions of $\hat x_i$, then the shifted-LP selection of Corollary~\ref{cor:inexact_stage1} satisfies
\begin{align*}
\sum_i f_i(\bar x_i)&\le f^*_{\mathrm{cvx}}+\Delta+\varepsilon_{\mathrm{obj}},\\
\Big\|\sum_iA_i\bar x_i-b\Big\|_\infty
&\le m_{\mathrm{eff}}\max_i\|A_i\|_\infty+\varepsilon_{\mathrm{feas}}.
\end{align*}
\end{corollary}
\begin{proof}
Corollary~\ref{cor:inexact_stage1} gives $\sum_i f_i(\bar x_i)\le G(\hat x)$. The envelope bound gives
$G(\hat x)\le\tilde G(\hat x)+\Delta$, while
$\tilde f^*_{\mathrm{cvx}}\le f^*_{\mathrm{cvx}}$ because $\tilde G\le G$ on $\mathcal F$. The residual
bound is~\eqref{eq:inexact_viol_bound}.
\end{proof}

\begin{remark}[Additive guarantee versus exact recovery]\label{rem:additive_vs_inclusion}
Corollary~\ref{cor:additive_value} asks nothing about $x^*$, and its hypothesis is met by construction on every Stage~2(a) path: full enumeration returns $\mathcal{V}_i=\mathcal{X}_i\cap\mathcal{N}(\tilde x_i^*)$, which contains the decomposition Lemma~\ref{lem:convex_decomposition} produces at $\tilde x_i^*$; an exactly solved column generation returns such a decomposition as its final basis; and so does the ellipsoid route of Remark~\ref{rem:poly_decomposition}. The $\Delta$ term can be removed, recovering the exact bound~\eqref{eq:main_obj_bound}, if the candidate sets also contain an LCE decomposition of the envelope minimizer $x^*$:
\begin{align}
\{v_{ij}\}_{j=1}^{k_i}\subseteq \mathcal{V}_i,\quad i=1,\ldots,N,\label{eq:vertex_inclusion}
\end{align}
with $v_{ij}$ the support of~\eqref{eq:decomposition_i} at $x^*$. Then $y=\alpha(x^*)$ is feasible with right-hand side $b$ and cost $f^*_{\mathrm{cvx}}$, and purifying it gives the selection. The residual is governed by the candidate diameters of Proposition~\ref{prop:violation_bound}. If $\mathcal V_i\subseteq\mathcal N(x_i^*)\cup\mathcal N(\tilde x_i^*)$ and $\|x_i^*-\tilde x_i^*\|_\infty<1$, then any two candidates satisfy $\|v-v'\|_\infty\le\|v-x_i^*\|_\infty+\|x_i^*-\tilde x_i^*\|_\infty+\|\tilde x_i^*-v'\|_\infty<3$, and being integer vectors they have $\|v-v'\|_\infty\le2$. Hence $D_i\le2$ and the residual bound at most doubles. The inclusion~\eqref{eq:vertex_inclusion} refers to $x^*$, however, and is certified by no path. 
A norm bound alone does not ensure that the candidate sets contain an LCE decomposition of $x^*$. Proposition~\ref{prop:proximity} is therefore an aggregate proximity statement, not evidence for the inclusion.
\end{remark}

\section{Numerical Illustration}\label{sec:numerical_illustration}

We run the two-stage algorithm on the inventory-placement problem of Example~\ref{ex:placement} with a weighted capacity coupling. 
The weighted coupling allows nonintegral relaxation polytopes, placing this formulation in the regime the method targets.
 The experiments check the two guarantees of Corollary~\ref{cor:eps_stage1} on every run, and measure the computational time of the method against three exact integer-programming benchmarks run under a common protocol: a piecewise-linear MILP, the same MILP with the coupling relaxed by our own residual, and a direct MIQP. A final set of runs takes the method alone to $N=2^{20}$ products.

\subsection{Problem Instance}\label{sec:num_setup}

A crossdock receives $q_i\in\mathbb{N}$ units of each of $N$ products during a day and must route them to $n$ downstream regions. Writing $x_i^j$ for the number of units of product $i$ sent to region $j$, conservation of each product's inflow imposes $\sum_j x_i^j=q_i$, so $\mathcal{X}_i=\{x\in\mathbb{N}^n:\sum_j x^j=q_i\}$ is the canonical domain of Section~\ref{sec:mconvexity}. Region $j$ holds on-hand inventory $y_i^j\ge 0$ of product $i$ and has mean demand $\mu_i^j>0$. The placement cost is the demand-weighted squared percentage error of the resulting inventory position,
\begin{align}
f_i(x_i)\;=\;\sum_{j=1}^n \mu_i^j\left(\frac{x_i^j+y_i^j-\mu_i^j}{\mu_i^j}\right)^{\!2}\;=\;\sum_{j=1}^n\frac{\big(x_i^j+y_i^j-\mu_i^j\big)^2}{\mu_i^j},\label{eq:placement_cost}
\end{align}
a separable convex quadratic restricted to $\mathcal{X}_i$, hence $M$-convex. Each unit of product $i$ consumes $w_i\in\mathbb{N}$ units of a shared arc resource, and the arc into region $j$ admits $b^j$ units, with $\sum_j b^j=\sum_i w_i q_i$. The coupling $\sum_i w_i x_i^j=b^j$ is problem~\eqref{eq:base_problem} with $m=n$ rows and $A_i=w_i I$. Every candidate difference satisfies $\mathbf 1^\top A_i(v-v')=0$, so $m_{\mathrm{eff}}\le m-1$ (Proposition~\ref{prop:violation_bound}), and the residual bound we test is $(n-1)\max_i w_i$.

The instance data $(q_i,\mu_i^j,y_i^j)$ are drawn from distributions fit to a shadow dataset, the weights $w_i$ uniformly from $\{1,\ldots,10\}$, and the capacities $b$ are an even split of $\sum_i w_iq_i$ (Appendix~\ref{app:data_generation}).
The comparison grid is $n\in\{4,8,16,32\}$ and $N\in\{2^8,\ldots,2^{14}\}$ with ten seeds per cell; the scaling runs add $N\in\{2^{15},2^{16},2^{17},2^{18}\}$ with ten seeds and $N\in\{2^{19},2^{20}\}$ with three.

\subsection{Implementation of the Two-Stage Algorithm}\label{sec:num_solvers}

Every step below operates on arrays indexed by agent, so a step over all $N$ agents is a handful of vectorized array operations rather than a loop; the one exception is purification, whose per-step work does not depend on $N$.

\paragraph{Stage 1: consensus ADMM on the envelopes.} Stage~1 is the two-block consensus ADMM of Appendix~\ref{app:admm_details} applied to the relaxation~\eqref{eq:convex_relaxation} with the envelopes $f_i^{**}$ themselves, not a surrogate. For the separable cost~\eqref{eq:placement_cost}, Proposition~\ref{rem:separable_decomposition} identifies $f_i^{**}$ with the sum over coordinates of the chord interpolants of $t\mapsto(t+y_i^j-\mu_i^j)^2/\mu_i^j$ at the integers, so each proximal step $\mathcal P_{i,\rho}$ separates once the multiplier $\nu_i$ of $\sum_jx_i^j=q_i$ is fixed. Each coordinate then has a closed-form minimizer, and $\nu_i$ is the root of a monotone piecewise-linear equation in one variable (Appendix~\ref{app:admm_impl}). The projection onto the coupling is the diagonal closed form of Appendix~\ref{app:admm_details}. Over-relaxation and residual balancing of the penalty $\rho$ \citep{boyd2011distributed} are used; the run stops when the primal and dual residuals meet the tolerances of Appendix~\ref{app:solver_settings}, which happened on every run.

\paragraph{Coupling projection.} ADMM stops on its consensus residuals, so its $x$-iterate $x^*$ satisfies the coupling only approximately: $r=\sum_iA_ix_i^*-b$ is small relative to $b$, but $\|r\|_\infty$ grows with $N$ at a fixed relative tolerance, and by~\eqref{eq:inexact_viol_bound} it adds directly to the residual bound. Before Stage~2 we therefore compute an approximate Euclidean projection $\tilde x$ of $x^*$ onto
\[
\Big\{x:\ \textstyle\sum_iA_ix_i=b,\ \ x_i\in\mathrm{conv}(\mathcal X_i),\ \ \lfloor x_i^*\rfloor\le x_i\le\lceil x_i^*\rceil\ \ \forall i\Big\},
\]
using a semismooth Newton method on the $m$-dimensional dual (Appendix~\ref{app:admm_impl}). The box keeps every agent in the integer cell of its Stage~1 point, so agents that are integral stay integral and no candidate set grows; if the boxes cannot absorb $r$, the projection drops them, which no run required. The Newton solve stops at a relative tolerance that leaves $\|\sum_iA_i\tilde x_i-b\|_\infty<1$. Set $\hat b:=\sum_iA_i\tilde x_i$. The decomposition weights reproduce $\hat b$, and purification preserves this aggregate. Corollary~\ref{cor:eps_stage1} therefore bounds the final coupling residual by $(n-1)\max_iw_i+\|\hat b-b\|_\infty<(n-1)\max_iw_i+1$. Since the final residual and $(n-1)\max_iw_i$ are integers, the returned point satisfies
\begin{align}
\Big\|\sum_iA_i\bar x_i-b\Big\|_\infty\le(n-1)\max_iw_i,\qquad \sum_if_i(\bar x_i)\le\sum_if_i^{**}(\tilde x_i).\label{eq:num_guarantees}
\end{align}
The value bound is measured against the computable envelope objective $\sum_if_i^{**}(\tilde x_i)$ at the approximately feasible point; optimality gaps are measured against the benchmarks below.

\paragraph{Stage 2(a): decompositions.} For each agent, Stage~2(a) enumerates $\mathcal X_i\cap\mathcal N(\tilde x_i)$ when it is small (Appendix~\ref{app:solver_settings}) and obtains the weights $\alpha_i$ from a per-agent LP; otherwise it runs column generation on the decomposition LP~\eqref{eq:decomposition_lp}, started from the offset sweep of Appendix~\ref{app:offset_sweep}. For the separable cost the offset sweep alone already returns a value-preserving decomposition (Proposition~\ref{rem:separable_decomposition}), and its basis is already optimal for column generation. The enumerate-or-price dispatch is the general $M$-convex implementation, kept so that the reported times include it; it is a small share of the run time. Most agents are integral after Stage~1, and candidate sets are small. The per-agent work is distributed over a pool of worker processes.

\paragraph{Stage 2(b): purification.} Stage~2(b) runs Algorithm~\ref{alg:purification} on $y=\alpha$ with right-hand side $\hat b=\sum_iA_i\tilde x_i$ and the upper bound $\hat m=m-1$ on $m_{\mathrm{eff}}$.
Each step works on a block of fractional agents just large enough that the per-agent simplex and coupling constraints restricted to the block have a nonzero null vector, takes that null vector as the direction, orients it so that the cost does not increase, and moves until a weight reaches zero. The block has $O(m)$ agents, so a step costs $O(m^3)$ independently of $N$. The walk stops when at most $m-1$ agents remain fractional, and each of those takes its cheapest support vertex, as in~\eqref{eq:selection_rule}.

\subsection{Benchmarks and Protocol}\label{sec:num_benchmarks}

\paragraph{Exact piecewise-linear MILP.} Each coordinate cost of~\eqref{eq:placement_cost} is replaced by its integer-point piecewise-linear interpolant, written with one continuous unit increment per integer step; integrality is imposed only on the $Nn$ aggregates $x_i^j$, and convexity makes the increments fill cheapest-first, so the model is an exact integer reformulation of~\eqref{eq:base_problem} (Appendix~\ref{app:benchmarks}).
Its LP relaxation equals the relaxation~\eqref{eq:convex_relaxation}, which is fractional here. We solve it with FICO Xpress~\citep{fico2024xpress}.

\paragraph{Relaxed-coupling MILP.} The two-stage point is allowed a coupling residual that the exact benchmark forbids, so its cost can fall slightly below the exact optimum, and the two do not solve the same problem. The relaxed benchmark is the same MILP with each coupling row replaced by $b^j-r\le\sum_iw_ix_i^j\le b^j+r$, where $r$ is the residual of our own solution on that seed. Our solution is feasible for it, so its certified lower bound is a lower bound on the cost of any solution with at most our residual, and the distance from our cost to that bound is an upper bound on our optimality gap within that residual.

\paragraph{SCIP MIQP.} The direct formulation: integer $x_i^j$, the coupling and the per-agent sums as linear equalities, and the convex quadratic cost as one epigraph constraint per term, solved with SCIP~\citep{bestuzheva2023scip}. This is what a modeler obtains without reformulating. We omit the plain MIQP in Xpress: its continuous relaxation is the QP relaxation, which lies below the envelope by up to $1/(4\mu_i^j)$ per fractional coordinate, and in our tests its gap stayed near $10^{-3}$ at the time limit even on the smallest cells, an order of magnitude above the stopping target.

\paragraph{Protocol.} All methods run on the same machine under the same rules (Appendix~\ref{app:solver_settings}). Each job is pinned to the same four physical cores ($8$ hardware threads). SCIP's branch and bound is sequential, so each SCIP job is pinned to one core. Benchmarks start cold, with no solution or information from the two-stage method. They stop when they certify a relative gap of $10^{-4}$ on the true cost, about the accuracy at which the two-stage method itself lands, or at $3600$\,s. Reported times are full wall-clock times: for the benchmarks, model construction plus solve; for the two-stage method, both stages, the projection, and all bookkeeping. A benchmark that fails to certify a cell within the limit is not run on any cell with larger or equal $n$ and $N$.

\subsection{Results}\label{sec:num_results}

\begin{figure}[htbp]
\centering
\includegraphics[width=0.8\textwidth]{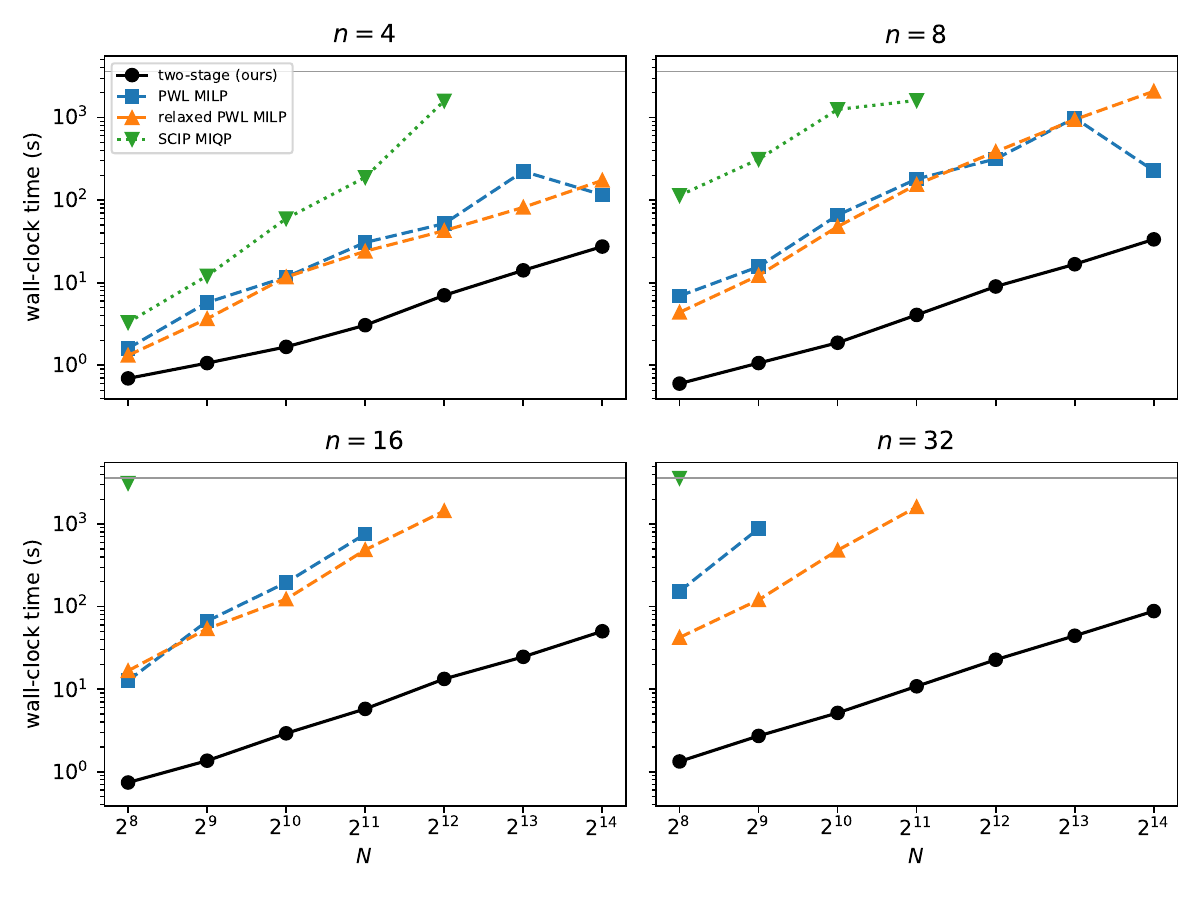}
\caption{Full wall-clock time against $N$ per region count $n$, for the two-stage method and the three exact benchmarks (medians over seeds, timed-out seeds counted at the $3600$\,s limit shown in gray). Benchmark curves stop where the benchmark was pruned.}\label{fig:time_vs_N}
\end{figure}

\begin{figure}[htbp]
\centering
\includegraphics[width=0.6\textwidth]{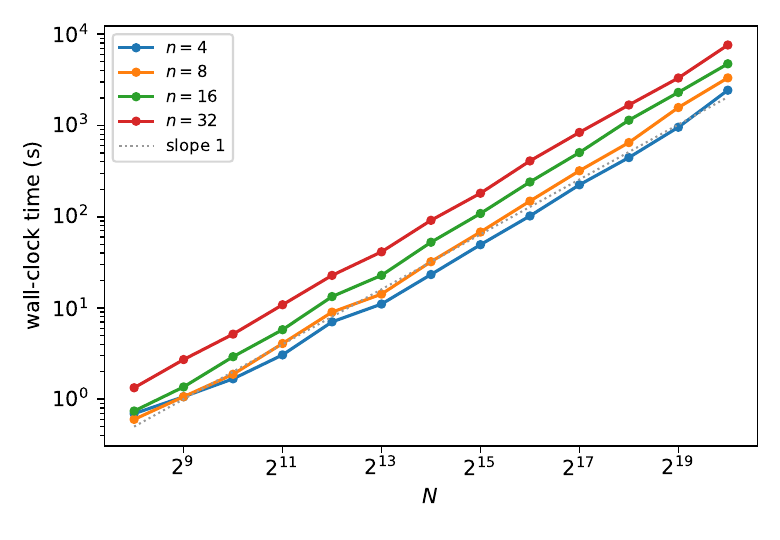}
\caption{Full wall-clock time of the two-stage method against $N$, one line per region count $n$ (medians over seeds), with a line of slope one for reference. Points with $N\le4{,}096$ are from the comparison grid of Table~\ref{tab:full_results}; points with $N\ge8{,}192$ are from the scaling runs of Table~\ref{tab:scale}.}\label{fig:ts_time_vs_N}
\end{figure}

Table~\ref{tab:full_results} and Figure~\ref{fig:time_vs_N} report the comparison grid, on which the two-stage method runs alongside the three benchmarks; Table~\ref{tab:scale} and Figure~\ref{fig:ts_time_vs_N} report the scaling runs of the method alone, up to $N=2^{20}$. The experiments have three aims: to check that the residual and value bounds~\eqref{eq:num_guarantees} hold on every run, and in particular that the residual does not grow with $N$; to measure the optimality gap of the integer solution against certified benchmarks; and to compare its running time with that of exact solvers. The results support the following observations.

\begin{table}[htbp]
\centering
\caption{Two-stage method against the three exact benchmarks, per $(n,N)$ over ten seeds. $r_{\max}$: largest coupling residual $\|\sum_iw_i\bar x_i-b\|_\infty$ over seeds, and the bound $(n-1)\max_iw_i$ of~\eqref{eq:num_guarantees}. ``time'': median full wall-clock time in seconds. ``gap'': mean over seeds of $(\sum_if_i(\bar x_i)-\underline f_r)/|\underline f_r|$, where $\underline f_r$ is the certified lower bound of the relaxed-coupling MILP, so the gap is an upper bound on the optimality gap within the residual $r$. Benchmark columns: median full time and, in parentheses, the number of seeds certified to a relative gap of $10^{-4}$ out of those run; a timed-out seed counts at the $3600$\,s limit, which understates the median, and is marked by $\dagger$; $\ddagger$ marks a SCIP cell where the solver aborted on some seed, which is excluded (Appendix~\ref{app:benchmarks}). ``pruned'': not run because a smaller cell timed out; ``--'': not run. The relaxed-coupling MILP was stopped before its largest cells for lack of memory (its model grows with $n\sum_iq_i$), so for $n=8$, $N=2^{14}$ and $n=32$, $N=2^{11}$ only the seeds shown were run, and cells beyond them are ``--'', not timeouts.}\label{tab:full_results}
\small\setlength{\tabcolsep}{4pt}
\begin{tabular}{rr rr rr rrr}
\toprule
& & \multicolumn{2}{c}{residual} & \multicolumn{2}{c}{two-stage} & \multicolumn{3}{c}{benchmark time (s), certified seeds} \\
\cmidrule(lr){3-4}\cmidrule(lr){5-6}\cmidrule(lr){7-9}
$n$ & $N$ & $r_{\max}$ & bound & time (s) & gap & MILP & relaxed MILP & SCIP MIQP \\
\midrule
4 & 256 & 16 & 30 & 0.7 & 0.021\% & 1.6 (10/10) & 1.3 (10/10) & 3.3 (10/10) \\
4 & 512 & 13 & 30 & 1.1 & 0.024\% & 5.8 (10/10) & 3.7 (10/10) & 12.1 (10/10) \\
4 & 1,024 & 15 & 30 & 1.7 & 0.046\% & 11.7 (10/10) & 11.7 (10/10) & 59.8 (10/10) \\
4 & 2,048 & 12 & 30 & 3.1 & 0.053\% & 30.7 (10/10) & 24.0 (10/10) & 189 (10/10) \\
4 & 4,096 & 13 & 30 & 7.0 & 0.082\% & 51.9 (10/10) & 42.6 (10/10) & 1,593$^\dagger$ (7/10) \\
4 & 8,192 & 18 & 30 & 14.1 & 0.066\% & 223 (10/10) & 81.7 (10/10) & pruned \\
4 & 16,384 & 17 & 30 & 27.4 & 0.079\% & 116 (10/10) & 173 (10/10) & pruned \\
\addlinespace
8 & 256 & 14 & 70 & 0.6 & 0.013\% & 6.9 (10/10) & 4.4 (10/10) & 114 (10/10) \\
8 & 512 & 15 & 70 & 1.1 & 0.011\% & 15.6 (10/10) & 12.2 (10/10) & 312 (10/10) \\
8 & 1,024 & 13 & 70 & 1.9 & 0.022\% & 65.7 (10/10) & 47.7 (10/10) & 1,253 (10/10) \\
8 & 2,048 & 16 & 70 & 4.1 & 0.023\% & 180 (10/10) & 154 (10/10) & 1,612$^\dagger$ (9/10) \\
8 & 4,096 & 18 & 70 & 9.0 & 0.037\% & 315 (10/10) & 386 (10/10) & pruned \\
8 & 8,192 & 18 & 70 & 16.7 & 0.035\% & 987 (10/10) & 945 (10/10) & pruned \\
8 & 16,384 & 19 & 70 & 33.5 & 0.052\% & 227 (10/10) & 2,080 (4/4) & pruned \\
\addlinespace
16 & 256 & 21 & 150 & 0.7 & 0.009\% & 12.7 (10/10) & 16.7 (10/10) & 3,127$^\dagger$ (5/10) \\
16 & 512 & 30 & 150 & 1.4 & 0.009\% & 66.7 (10/10) & 53.9 (10/10) & pruned \\
16 & 1,024 & 26 & 150 & 2.9 & 0.015\% & 196 (10/10) & 123 (10/10) & pruned \\
16 & 2,048 & 23 & 150 & 5.8 & 0.015\% & 757$^\dagger$ (9/10) & 486 (10/10) & pruned \\
16 & 4,096 & 24 & 150 & 13.3 & 0.024\% & pruned & 1,441 (10/10) & pruned \\
16 & 8,192 & 23 & 150 & 24.6 & -- & pruned & -- & pruned \\
16 & 16,384 & 23 & 150 & 50.3 & -- & pruned & -- & pruned \\
\addlinespace
32 & 256 & 26 & 310 & 1.3 & 0.011\% & 153 (10/10) & 42.3 (10/10) & 3,600$^\dagger$ (0/1) \\
32 & 512 & 21 & 310 & 2.7 & 0.007\% & 878$^\dagger$ (8/10) & 121 (10/10) & pruned \\
32 & 1,024 & 28 & 310 & 5.2 & 0.008\% & pruned & 481 (10/10) & pruned \\
32 & 2,048 & 28 & 310 & 10.8 & 0.007\% & pruned & 1,610 (8/8) & pruned \\
32 & 4,096 & 29 & 310 & 22.7 & -- & pruned & -- & pruned \\
32 & 8,192 & 28 & 310 & 44.4 & -- & pruned & -- & pruned \\
32 & 16,384 & 28 & 310 & 88.2 & -- & pruned & -- & pruned \\
\bottomrule
\end{tabular}
\end{table}

\paragraph{The guarantees hold on every run.} Over the $280$ comparison runs and the $264$ scaling runs, the coupling residual stays below $0.60$ of the bound $(n-1)\max_iw_i$ and does not grow with $N$: from $N=2^8$ to $2^{20}$ its maximum is $11$ to $18$ at $n=4$, against a bound of $30$, and $20$ to $29$ at $n=32$, against $310$. The per-product constraints $\sum_j\bar x_i^j=q_i$ hold exactly, and the cost lies below $\sum_if_i^{**}(\tilde x_i)$ on every run, by $10^{-6}$ to $9\times10^{-5}$ in relative terms.

\paragraph{The solution is near-optimal.} For cells with certified bounds from the relaxed-coupling MILP, the mean relative gap reported in Table~\ref{tab:full_results} ranges from $0.007\%$ to $0.082\%$; the gap is largest at $n=4$ and below $0.02\%$ for $n\ge16$. Against the exact-coupling MILP the differences are of the same size, and where our cost is the lower one (some seeds at $n\ge16$), the relaxed benchmark shows that this comes from the residual, not from a better solution.

\paragraph{The method is much faster than every benchmark, increasingly so with $n$ and $N$.} The speedup over the exact MILP grows from $2$ to $16$ at $n=4$ to $115$ to $323$ at $n=32$, and the exact MILP fails to certify within the hour from $N=2^{11}$ at $n=16$ and from $N=2^9$ at $n=32$, cells the two-stage method covers in under $90$\,s. Relaxing the coupling makes the MILP three to seven times faster at $n=32$, but our method remains $20$ to $150$ times faster for $n\ge16$. SCIP is the slowest method: $5$ to $60$ times slower than ours at $n=4$ and $190$ to $670$ times at $n=8$, it times out from $N=2^{12}$ at $n=4$ and $N=2^{11}$ at $n=8$, and for $n\ge16$ it does not certify every seed even at $N=2^8$. Nor do the MILPs find a solution as good as ours much before they certify: their first such incumbent usually appears close to the certification time, and even their first feasible incumbent arrives after the complete two-stage run. All but five of the certified MILP solves close at the root node, so their time is set by when the root heuristics find a near-optimal incumbent rather than by branching; this is why the MILP times are not monotone in $N$ (at $n=8$, the exact MILP is faster at $N=2^{14}$ than at $N=2^{13}$).

\paragraph{The running time grows almost linearly in $N$.} The time grows by a factor of $2.0$ to $2.5$ per doubling of $N$ and by three to four from $n=4$ to $n=32$; at $N=2^{20}$ a run takes $40$ minutes at $n=4$ and just over two hours at $n=32$, $70\%$ to $86\%$ of it in Stage~1 and under $2\%$ in Stage~2(a), which takes at most $5\%$ on every comparison cell with $N\ge 4096$.

\begin{table}[htbp]
\centering
\caption{Two-stage method alone at large $N$: number of seeds, largest coupling residual and its bound, and median times of Stage~1 (including the coupling projection), Stage~2, and the whole run, in seconds.}\label{tab:scale}
\small
\begin{tabular}{rr rrr rrr}
\toprule
$n$ & $N$ & seeds & $r_{\max}$ & bound & Stage~1 (s) & Stage~2 (s) & total (s) \\
\midrule
4 & 8,192 & 10 & 18 & 30 & 8.3 & 2.7 & 11.0 \\
4 & 16,384 & 10 & 17 & 30 & 18.1 & 5.2 & 23.2 \\
4 & 32,768 & 10 & 12 & 30 & 39.4 & 9.5 & 49.3 \\
4 & 65,536 & 10 & 12 & 30 & 82.8 & 18.2 & 102 \\
4 & 131,072 & 10 & 16 & 30 & 185 & 37.6 & 224 \\
4 & 262,144 & 10 & 12 & 30 & 366 & 76.5 & 445 \\
4 & 524,288 & 3 & 16 & 30 & 786 & 164 & 957 \\
4 & 1,048,576 & 3 & 11 & 30 & 2,048 & 381 & 2,438 \\
\addlinespace
8 & 8,192 & 10 & 18 & 70 & 10.1 & 4.0 & 14.2 \\
8 & 16,384 & 10 & 19 & 70 & 24.3 & 7.5 & 32.1 \\
8 & 32,768 & 10 & 18 & 70 & 53.8 & 13.5 & 67.9 \\
8 & 65,536 & 10 & 14 & 70 & 121 & 25.9 & 148 \\
8 & 131,072 & 10 & 15 & 70 & 263 & 54.7 & 319 \\
8 & 262,144 & 10 & 15 & 70 & 544 & 110 & 650 \\
8 & 524,288 & 3 & 13 & 70 & 1,353 & 246 & 1,574 \\
8 & 1,048,576 & 3 & 15 & 70 & 2,718 & 589 & 3,326 \\
\addlinespace
16 & 8,192 & 10 & 23 & 150 & 16.1 & 6.3 & 22.8 \\
16 & 16,384 & 10 & 23 & 150 & 40.6 & 11.5 & 52.5 \\
16 & 32,768 & 10 & 23 & 150 & 86.1 & 22.2 & 108 \\
16 & 65,536 & 10 & 29 & 150 & 200 & 39.7 & 240 \\
16 & 131,072 & 10 & 27 & 150 & 417 & 85.0 & 505 \\
16 & 262,144 & 10 & 23 & 150 & 973 & 171 & 1,143 \\
16 & 524,288 & 3 & 17 & 150 & 1,932 & 362 & 2,306 \\
16 & 1,048,576 & 3 & 14 & 150 & 3,881 & 854 & 4,744 \\
\addlinespace
32 & 8,192 & 10 & 28 & 310 & 29.0 & 12.1 & 41.2 \\
32 & 16,384 & 10 & 28 & 310 & 69.4 & 21.6 & 91.2 \\
32 & 32,768 & 10 & 24 & 310 & 141 & 39.3 & 181 \\
32 & 65,536 & 10 & 21 & 310 & 336 & 74.7 & 409 \\
32 & 131,072 & 10 & 27 & 310 & 686 & 159 & 839 \\
32 & 262,144 & 10 & 26 & 310 & 1,363 & 310 & 1,679 \\
32 & 524,288 & 3 & 21 & 310 & 2,738 & 638 & 3,317 \\
32 & 1,048,576 & 3 & 20 & 310 & 6,389 & 1,504 & 7,628 \\
\bottomrule
\end{tabular}
\end{table}

\section{Conclusion}\label{sec:conclusion}

We studied the minimization of a sum of $N$ $M$-convex functions coupled by $m$ linear constraints, a
problem for which no polynomial-time algorithm is known once $N\ge3$. Replacing each function by its lower
convex envelope gives a separable convex relaxation, and the envelope's agreement with the function at
integer points lets an optimal relaxed solution be rounded at no cost in objective: local
Carath\'eodory decompositions express each agent's point through nearby integer points, and an LP
purification leaves at most $m$ agents fractional. The resulting integer solution costs no more than the
relaxation value and violates the coupling by at most $m\max_i\|A_i\|_\infty$, independently of $N$. Both
stages can be carried out with polynomially many value-oracle calls, the relaxation can be solved in
parallel across agents through a proximal oracle, and a strongly convex surrogate can replace the envelope
in Stage~1 at an additive cost controlled by the envelope gap.

On an inventory placement problem with up to about $10^6$ products, the coupling residual stayed below
$0.6$ of the bound and did not grow with $N$. Where certified benchmark bounds were available, the mean
relative gaps were at most about $0.08\%$. The method was faster than the integer-programming benchmarks
by a factor that grew with the number of regions. The relaxation dominates the running time, which grows
almost linearly in $N$.

Several questions remain open. The residual bound is worst-case, and refining it under the structural
assumptions of \citet{bi2020refined} and \citet{dey2026asymptotically} could give tighter guarantees for
smooth or block-structured instances. When the coupling must hold exactly, the residual has to be absorbed
by a repair step or by slack in the constraints, and we have not analyzed the cost of doing so. Our
experiments use a separable cost, for which the decomposition step is trivial; nonseparable $M$-convex
costs, such as those arising from substitutable demand, would exercise the pricing oracle. Finally, the
envelope gap of a strongly convex surrogate grows with $N$, so a surrogate whose loss stays bounded
independently of $N$, while still allowing linear-rate Stage~1 solvers, would remove the main drawback of
that variant.

\bibliographystyle{plainnat}
\bibliography{bibliography}

\appendix

\section{Consensus ADMM Subproblem Derivation}\label{app:admm_details}

For the LCE relaxation, Lemma~\ref{lem:lce_saddle} supplies $\lambda^*$ with
$A_i^\top\lambda^*\in\partial f_i^{**}(x_i^*)$ at every optimal $x^*$, so $(x,z,\gamma)=(x^*,x^*,A^\top\lambda^*)$
is a saddle point of the consensus Lagrangian below. For a surrogate
relaxation, assume that a saddle point exists. Throughout this section, the proximal subproblems are
solved exactly. The displayed rate is for standard ADMM, without over-relaxation.

The 2-block consensus ADMM for~\eqref{eq:convex_relaxation} introduces a copy $z=(z_1,\ldots,z_N)\in\mathbb{R}^{Nn}$ and absorbs the coupling constraint $Az=b$ (with $A=[A_1\;\cdots\;A_N]$) into the indicator $\iota_{\mathcal{C}}$ on $\mathcal{C}=\{z:Az=b\}$. The scaled-dual augmented Lagrangian
\begin{align*}
\mathcal{L}_\rho(x,z,u)=\sum_{i=1}^N f_i^{**}(x_i)+\iota_{\mathcal{C}}(z)
+\frac{\rho}{2}\|x-z+u\|_2^2-\frac{\rho}{2}\|u\|_2^2
\end{align*}
yields
\begin{align}
x_i^{k+1}&=\mathcal{P}_{i,\rho}(z_i^k-u_i^k),\quad i=1,\ldots,N,\label{eq:cadmm_x_step}\\
z^{k+1}&=\Pi_{\mathcal{C}}(x^{k+1}+u^k),\label{eq:cadmm_z_step}\\
u^{k+1}&=u^k+x^{k+1}-z^{k+1}.\label{eq:cadmm_u_step}
\end{align}
The $N$ proximal calls run in parallel. With a saddle point and exact proximal solves, standard
2-block ADMM with any fixed $\rho>0$ drives the residual $x^k-z^k$ to zero and the objective to its optimal
value \citep{boyd2011distributed}; convergence of the iterates themselves follows from the
Douglas--Rachford interpretation of ADMM \citep[Thm.~8]{eckstein1992douglas}.

\paragraph{$x$-update.} Equation~\eqref{eq:cadmm_x_step} is exactly the proximal-oracle call
$\mathcal{P}_{i,\rho}(z_i^k-u_i^k)$ from~\eqref{eq:prox_oracle}, independently for each agent. We treat this
map as primitive here and in the main presentation, so the iteration bound below is in units of exact
proximal calls ($N$ per iteration).
Remark~\ref{rem:prox_value_oracle} sketches how to implement each call with polynomially many evaluations
of the underlying $M$-convex function $f_i$.

\paragraph{$z$-update.} The projection onto $\mathcal{C}$ has closed form
\begin{align}
\Pi_{\mathcal{C}}(g)=g-A^\top(AA^\top)^\dagger(Ag-b),\label{eq:cadmm_proj_general}
\end{align}
where $\dagger$ denotes the Moore-Penrose pseudoinverse. A rank-revealing factorization of $AA^\top$ is
computed once; each iteration then uses two matvecs with $A$ and one solve on its range.

When each $A_i=\mathrm{diag}(a_i)$, $AA^\top$ is diagonal with $(AA^\top)_{jj}=\sum_i(a_i^j)^2$, and~\eqref{eq:cadmm_proj_general} factors per-coordinate:
\begin{align*}
z_i^{k+1,j}=g_i^{k,j}-\frac{a_i^j}{\sum_{i'}(a_{i'}^j)^2}\Big(\sum_{i'}a_{i'}^j g_{i'}^{k,j}-b^j\Big),\quad i=1,\ldots,N,\;j=1,\ldots,n,
\end{align*}
with $g^k:=x^{k+1}+u^k$ for every nonzero row of $A$; zero rows are redundant on a feasible instance and
are removed. For $A_i=I$ this simplifies to
$z_i^{k+1,j}=g_i^{k,j}-\tfrac{1}{N}(\sum_{i'} g_{i'}^{k,j}-b^j)$.

\paragraph{Dual update.} $u^{k+1}=u^k+x^{k+1}-z^{k+1}$. The unscaled multiplier on $x=z$ in the sign
convention used below is $\gamma^k=-\rho u^k$. For $k\ge1$, define
\begin{align}
\lambda^k:=-\rho(AA^\top)^\dagger Au^k.\label{eq:finite_multiplier}
\end{align}
Indeed, $u^{k+1}=g^k-z^{k+1}$ is a projection residual and therefore lies in
$\mathrm{range}(A^\top)$. Hence $A^\top\lambda^k=-\rho u^k=\gamma^k$. If ADMM converges,
$\lambda^k\to\lambda^*$ and optimality of the limiting $x$-step gives
$A_i^\top\lambda^*\in\partial f_i^{**}(x_i^*)$.
This is the Stage~1 multiplier in the sign convention of~\eqref{eq:lagrangian}. If $A$ is rank deficient,
\eqref{eq:finite_multiplier} selects the minimum-norm representative. Adding an element of
$\ker A^\top$ changes neither $A^\top\lambda^*$ nor $\lambda^{*\top}b$, because feasibility implies
$b\in\mathrm{range}(A)$.

\paragraph{Ergodic output.} Define
\begin{align}
\bar x^K:=\frac1K\sum_{k=1}^Kx^k,\qquad
\bar z^K:=\frac1K\sum_{k=1}^Kz^k,\label{eq:ergodic_averages}
\end{align}
and $F(x):=\sum_if_i^{**}(x_i)$. The ergodic variational-inequality estimate of
\citet[Thm.~4.1]{he2012convergence} yields
\begin{align}
F(\bar x^K)-F(x^*)\le \frac{C_{\mathrm{ADMM}}}{K},\qquad
\|A\bar x^K-b\|_\infty\le\frac{C_{\mathrm{ADMM}}}{K},
\label{eq:admm_ergodic_bounds}
\end{align}
for a constant $C_{\mathrm{ADMM}}$ that depends on $\rho$, on the initialization $(z^0,u^0)$ and on an
optimal primal-dual pair; it is an analytic complexity parameter, not a stopping certificate.

\paragraph{Over-relaxation (optional).} Replacing $x^{k+1}$ by $\hat x^{k+1}:=\alpha\,x^{k+1}+(1-\alpha)\,z^k$ in both the $z$- and $u$-updates yields the over-relaxed variant; $\alpha\in[1.5,1.8]$ typically gives empirical speedups~\citep{boyd2011distributed}.

\section{Value-Oracle Implementations of Stage~1}\label{app:stage1_complexity}

This section fixes the rational input model and supplies the exact Stage~1
implementation cited in Section~\ref{sec:stage1}: Theorem~\ref{thm:exact_stage1}
solves the convex relaxation exactly through an exponentially large
configuration LP with polynomially many value-oracle calls.
Section~\ref{app:prox_value_oracle} then records, as a side result with a
proof sketch, that the proximal primitive used by distributed consensus ADMM
admits feasible approximations through polynomially many evaluations of the
underlying $M$-convex function.

We use the following input model. Write $\langle z\rangle$ for the binary encoding length of a rational
object $z$. Each value oracle is defined on all of $\mathbb{Z}^n$, returns $+\infty$ outside the nonempty
bounded $M$-convex set $\mathcal{X}_i=\mathrm{dom}\,f_i$, and returns rational finite values on
$\mathcal X_i$. For each agent, a point $v_i^0\in\mathcal X_i$ and a positive common denominator $D_i$
for its finite values are supplied. The matrices $A_i$, vector $b$, and algorithmic parameters are
rational. The following promise bounds are also supplied:
\begin{align}
\|v\|_\infty\le R_{\mathrm{dom}},\qquad
|f_i(v)|\le W_{\mathrm{glob}}
\quad\text{for every }i\text{ and }v\in\mathcal X_i.\label{eq:global_bounds}
\end{align}
Let $L_{\mathrm{data}}$ be the total encoding length of the rational data, the $v_i^0$, the $D_i$, the
promise bounds, and any fixed rational algorithmic parameters. All bounds below count value-oracle calls
and arithmetic operations on rationals of polynomial encoding length. When each $f_i$ can be evaluated
exactly in time polynomial in the encoding length of its parameters and of the query point, as for the
placement cost~\eqref{eq:placement_cost}, they are polynomial running-time bounds in the ordinary sense.

A linear perturbation of an $M$-convex function remains $M$-convex. After clearing the denominators of
the function values and the perturbation, full-domain minimization is polynomial by the general scaling
algorithm of \citet{tamura2005coordinatewise}; see also the related scaling approach of \citet{moriguchi2002scaling}. Its value-oracle count
is polynomial in $n$, $\log(1+R_{\mathrm{dom}})$, $\log(1+W_{\mathrm{glob}})$, and the encoding length
of the perturbation and cleared denominators. Indeed,
$\max_{u,v\in\mathcal X_i}\|u-v\|_\infty\le2R_{\mathrm{dom}}$. Scaling is needed because, unlike the
unit-neighborhood pricing problem in Theorem~\ref{thm:poly_extension}, the full domain need not have
unit diameter.

The same full-domain oracle supplies the initialization required by Theorem~\ref{thm:poly_extension}.
Given rational $x\in\mathrm{conv}(\mathcal{X}_i)$, minimize over $v\in\mathcal{X}_i$ the separable convex
box-distance
\begin{align}
d_x(v):=\sum_{j=1}^n\left[\max\{\lfloor x^j\rfloor-v^j,0\}
+\max\{v^j-\lceil x^j\rceil,0\}\right].\label{eq:local_phase_one}
\end{align}
Adding a separable convex function to the indicator of an $M$-convex set preserves $M$-convexity, so the
scaling algorithm applies from $v_i^0$ \citep[Ch.~6]{murota2003discrete}. Integral convexity implies
$\mathcal{X}_i\cap\mathcal{N}(x)\ne\varnothing$, hence the minimum is zero and any minimizer is a valid
local starting point.

\subsection{Exact solution through the configuration LP}\label{app:configuration_lp}

Introduce one variable $y_{iv}$ for each $v\in\mathcal{X}_i$ and consider
\begin{align}
(\mathrm{CLP}):\quad
\min_y\quad &\sum_{i=1}^N\sum_{v\in\mathcal{X}_i}f_i(v)y_{iv}\label{eq:configuration_lp}\\
\mathrm{s.t.}\quad
&\sum_{v\in\mathcal{X}_i}y_{iv}=1,\quad i=1,\ldots,N,\nonumber\\
&\sum_{i=1}^N\sum_{v\in\mathcal{X}_i}A_iv\,y_{iv}=b,\qquad y\ge0.\nonumber
\end{align}

\begin{lemma}[Configuration-LP equivalence]\label{lem:configuration_equiv}
Problem~\eqref{eq:configuration_lp} and the convex relaxation~\eqref{eq:convex_relaxation} have the same
optimal value. Moreover, an optimal $y^*$ induces an optimal point
$x_i^*=\sum_v y_{iv}^*v$ of~\eqref{eq:convex_relaxation}, with
\begin{align}
\sum_{v\in\mathcal{X}_i}f_i(v)y_{iv}^*=f_i^{**}(x_i^*)\qquad\text{for every }i.\label{eq:config_agent_identity}
\end{align}
\end{lemma}

\begin{proof}
Any feasible $y$ induces $x_i=\sum_vy_{iv}v\in\mathrm{conv}(\mathcal{X}_i)$ satisfying the coupling
constraint. By the definition of the lower convex envelope,
$f_i^{**}(x_i)\le\sum_vf_i(v)y_{iv}$, so the configuration-LP value is at least
$f^*_{\mathrm{cvx}}$. Conversely, any feasible $x$ of~\eqref{eq:convex_relaxation} admits, independently
for each agent, a finite convex decomposition satisfying
$x_i=\sum_vy_{iv}v$ and $f_i^{**}(x_i)=\sum_vf_i(v)y_{iv}$. These decompositions form a feasible point
of~\eqref{eq:configuration_lp} of the same value. The two optima are therefore equal. At an optimal $y^*$,
the nonnegative gaps
$\sum_vf_i(v)y_{iv}^*-f_i^{**}(x_i^*)$ sum to zero, which proves~\eqref{eq:config_agent_identity}.
\end{proof}

\begin{lemma}[Multipliers for the LCE relaxation]\label{lem:lce_saddle}
The LP dual of~\eqref{eq:configuration_lp} is
\begin{align}
\max_{\eta\in\mathbb{R}^N,\,\lambda\in\mathbb{R}^m}\quad
&\sum_{i=1}^N\eta_i+\lambda^\top b\label{eq:configuration_dual}\\
\mathrm{s.t.}\quad
&\eta_i+\lambda^\top A_iv\le f_i(v),
\qquad v\in\mathcal{X}_i,\quad i=1,\ldots,N.\nonumber
\end{align}
If~\eqref{eq:convex_relaxation} is feasible, let $(\eta^*,\lambda^*)$ be optimal for~\eqref{eq:configuration_dual}
and let $x^*$ be any optimal point of~\eqref{eq:convex_relaxation}. Then
\begin{align}
A_i^\top\lambda^*\in\partial f_i^{**}(x_i^*),\qquad i=1,\ldots,N,\label{eq:configuration_kkt}
\end{align}
so each $x_i^*$ minimizes $f_i^{**}(x_i)-\lambda^{*\top}A_ix_i$, and $(x^*,\lambda^*)$ is an optimal
primal-dual pair in the sign convention of~\eqref{eq:lagrangian}.
\end{lemma}

\begin{proof}
Dual feasibility says that the affine function $x\mapsto\eta_i^*+\lambda^{*\top}A_ix$ lies below $f_i$ on
$\mathcal X_i$. It is therefore an affine minorant of $\bar f_i$, and hence
\begin{align}
\eta_i^*+\lambda^{*\top}A_ix_i\le f_i^{**}(x_i)\qquad(x_i\in\mathbb R^n).
\label{eq:configuration_supporting_plane}
\end{align}
Summing at $x^*$ and using $\sum_iA_ix_i^*=b$ gives
$\sum_i(\eta_i^*+\lambda^{*\top}A_ix_i^*)=\sum_i\eta_i^*+\lambda^{*\top}b$. By LP strong duality and
Lemma~\ref{lem:configuration_equiv}, this equals $f^*_{\mathrm{cvx}}=\sum_if_i^{**}(x_i^*)$. The
nonnegative gaps $f_i^{**}(x_i^*)-\eta_i^*-\lambda^{*\top}A_ix_i^*$ therefore sum to zero, so each
vanishes. Each affine minorant in~\eqref{eq:configuration_supporting_plane} is thus tight at $x_i^*$,
which is~\eqref{eq:configuration_kkt}.
\end{proof}

\begin{theorem}[Exact value-oracle-polynomial Stage~1]\label{thm:exact_stage1}
Under the input model above, a feasible instance of~\eqref{eq:convex_relaxation} can be solved exactly with a
number of evaluations of the $f_i$ and rational arithmetic operations polynomial in $N$, $n$, $m$,
$\log(1+R_{\mathrm{dom}})$, $\log(1+W_{\mathrm{glob}})$ and $L_{\mathrm{data}}$. The
algorithm returns a rational optimal primal-dual pair $(x^*,\lambda^*)$ of polynomial encoding length.
\end{theorem}

\begin{proof}
We use the separation--optimization and oracle-dual recovery machinery of
Theorem~\ref{thm:poly_extension} and Remark~\ref{rem:poly_decomposition}.
By Lemma~\ref{lem:configuration_equiv}, feasibility of the convex relaxation
makes~\eqref{eq:configuration_lp} feasible. Since every bounded integer
domain $\mathcal X_i$ is finite, this is a finite LP with a finite attained
optimum. Its dual is~\eqref{eq:configuration_dual}, and LP strong duality gives dual attainment and
equality of the optimal values.

Given rational $(\eta,\lambda)$, separation checks, for each agent, whether
\begin{align}
\min_{v\in\mathcal{X}_i}\{f_i(v)-\lambda^\top A_iv\}\ge\eta_i.
\label{eq:configuration_pricing}
\end{align}
The perturbation $A_i^\top\lambda$ has encoding length polynomial in
$L_{\mathrm{data}}+\langle\lambda\rangle$. After clearing its denominators
together with $D_i$, the full-domain scaling algorithm described above,
initialized at $v_i^0$, uses polynomially many evaluations of $f_i$ and
rational operations in
\[
n,\quad \log(1+R_{\mathrm{dom}}),\quad
\log(1+W_{\mathrm{glob}}),\quad
L_{\mathrm{data}},\quad\text{and}\quad
\langle\lambda\rangle .
\]
Any trial point satisfying $\|v\|_\infty>R_{\mathrm{dom}}$ can be assigned
value $+\infty$ without calling the oracle. Hence every actual value-oracle
query has polynomial encoding length. A violating minimizer supplies the
inequality with normal $(e_i,A_iv)$; otherwise all constraints for agent $i$
hold. Thus one separation call uses at most $N$ pricing calls.

For $v\in\mathcal X_i$, the normal $(e_i,A_iv)$, right-hand side $f_i(v)$,
and dual objective $(\mathbf 1,b)$ have encoding length bounded by a polynomial $\varphi$ in $n$, $m$,
$\log(1+R_{\mathrm{dom}})$, $\log(1+W_{\mathrm{glob}})$ and $L_{\mathrm{data}}$, independent of $v$.
Hence~\eqref{eq:configuration_dual} is a well-described rational polyhedron of facet complexity at most
$\varphi$, with a strong separation oracle whose cost is polynomial in the input and query encoding
lengths. It is full-dimensional, since decreasing every $\eta_i$ from a feasible point makes all
constraints strict. It may be unbounded and may contain a lineality space (for $A_i=I$ on fixed-sum
domains, $\lambda\mapsto\lambda+c\mathbf 1$ together with $\eta_i\mapsto\eta_i-cq_i$ preserves every constraint and,
since feasibility of~\eqref{eq:configuration_lp} forces $\mathbf 1^\top b=\sum_iq_i$, also the objective), but its objective is bounded above because the primal~\eqref{eq:configuration_lp} is
feasible. The strong-separation/strong-optimization theorem for well-described polyhedra
\citep[Thm.~6.4.9]{grotschel1988geometric} therefore returns an exact
optimal full-dual vector $(\eta^*,\lambda^*)$ of polynomial encoding length
using polynomially many separation calls.

The constructive oracle-dual recovery result
\citep[Lem.~6.5.15]{grotschel1988geometric} makes polynomially many additional
separation calls and produces a polynomial-size set of generated inequalities
such that the configuration LP restricted to the corresponding columns
attains the full optimum. The inequality indexed by $(i,v)$ corresponds to
the column $(e_i,A_iv)$, while the primal right-hand side is
$(\mathbf 1,b)$. Solve this restricted LP exactly and choose a basic optimum
$y^*$. It is feasible for the full configuration LP, attains its optimum,
and has at most $N+m$ positive entries. Since the restricted LP has
polynomial dimensions and polynomial-size rational data, determinant bounds
give polynomial encoding length for $y^*$, and hence for
\[
x_i^*:=\sum_v y_{iv}^*v.
\]
We retain the full-dual optimizer $(\eta^*,\lambda^*)$, since an arbitrary
optimizer of the restricted dual need not satisfy the omitted inequalities.

Lemma~\ref{lem:configuration_equiv} makes $x^*$ optimal for the convex
relaxation. Since $(\eta^*,\lambda^*)$ is optimal for the full
dual~\eqref{eq:configuration_dual}, Lemma~\ref{lem:lce_saddle} gives~\eqref{eq:configuration_kkt},
so $(x^*,\lambda^*)$ is an optimal primal-dual pair in the sign convention
of~\eqref{eq:lagrangian}.

The optimization and recovery procedures make polynomially many separation
calls, each using at most $N$ pricing calls of polynomial value-oracle and
rational-arithmetic complexity. The exact restricted-LP solve is also
polynomial in its dimensions and encoding length. This proves the stated
value-oracle and rational-operation bounds.
\end{proof}

\begin{remark}[From exact Stage~1 to the full algorithm]\label{rem:exact_end_to_end}
The configuration-LP support is not used for reconstruction: its points may be far apart, which would
break the unit-diameter property behind Theorem~\ref{thm:main}. Instead, for each rational $x_i^*$
returned by Theorem~\ref{thm:exact_stage1}, Theorem~\ref{thm:poly_extension} and
Remark~\ref{rem:poly_decomposition} recover a value-preserving decomposition in
$\mathcal{X}_i\cap\mathcal{N}(x_i^*)$, initialized as in Section~\ref{sec:mconvexity} for the canonical
domain and by~\eqref{eq:local_phase_one} in general, with local bound
$\bar W_i:=\lceil D_iW_{\mathrm{glob}}\rceil$. The selection LP, with at most $Nn$ variables and $N+m$
equality constraints, is solved exactly by rational linear programming, returning an optimal basic
solution. Every step of the two-stage method therefore has value-oracle-polynomial complexity, and
Theorem~\ref{thm:main} applies with its stated residual.
\end{remark}

\subsection{Approximating the proximal oracle from values of \texorpdfstring{$f_i$}{f\_i}}\label{app:prox_value_oracle}

Fix an agent and suppress its index. For rational $\rho>0$ and $v\in\mathbb Q^n$, with
$R_v(x):=\frac{\rho}{2}\|x-v\|_2^2$, which is finite on
$\mathbb R^n$, and $(f^{**})^*=f^*$, Fenchel duality \citep[Thm.~31.1]{rockafellar1970convex} gives
\begin{align}
\min_x\{f^{**}(x)+R_v(x)\}=\frac{\rho}{2}\|v\|_2^2-\min_{p\in\mathbb R^n}\Psi(p),\qquad
\Psi(p):=f^*(p)+\frac{1}{2\rho}\|p-\rho v\|_2^2,\label{eq:prox_dual}
\end{align}
and $\mathcal P_\rho(v)=v-p^*/\rho$ at the unique minimizer $p^*$ of the $1/\rho$-strongly convex
$\Psi$. Only values of $f$ enter: $f^*(p)=-\min_{w\in\mathcal X}\{f(w)-p^\top w\}$ is one full-domain
$M$-convex minimization, solved by the scaling algorithm above after clearing denominators with $D_i$,
whose minimizer $w(p)$ also gives the subgradient $w(p)+(p-\rho v)/\rho$ of
$\Psi$.

\begin{remark}[Value-oracle implementation of the prox]\label{rem:prox_value_oracle}
Fix rational $\varepsilon\in(0,1]$. One can compute a rational $\hat x\in\mathrm{conv}(\mathcal X)$ with
$f^{**}(\hat x)+R_v(\hat x)\le\min_x\{f^{**}(x)+R_v(x)\}+\varepsilon$, using a number of evaluations of
$f$ and rational arithmetic operations polynomial in $n$, $\log(1/\varepsilon)$, and the encoding lengths
of $R_{\mathrm{dom}}$, $W_{\mathrm{glob}}$, the agent's value denominator $D_i$, $\rho$, and $v$. By $\rho$-strong convexity this gives
$\|\hat x-\mathcal P_\rho(v)\|_2\le\sqrt{2\varepsilon/\rho}$.

We sketch the argument. Set $\bar R:=\max\{1,R_{\mathrm{dom}}\}$. The minimizer
$p^*=\rho(v-\mathcal P_\rho(v))$ lies in the ball $\|p-\rho v\|_2\le2\rho\sqrt n\bar R$, at distance at
least $\rho\sqrt n\bar R$ from its boundary, and the subgradients above are bounded by $3\sqrt n\bar R$ on
this ball. A rounded central-cut ellipsoid run of fixed length $O(n^2\log(2+\rho n\bar R^2/\varepsilon))$,
with centers of polynomial encoding length and stopped early if a subgradient vanishes, therefore produces
an $\varepsilon$-optimal queried center \citep[Chs.~3--4]{grotschel1988geometric}.

Let $S\subseteq\mathcal X$ collect the minimizers $w(p)$ returned during the run, and let $\Psi_S$ be
$\Psi$ with $f^*$ restricted to $S$. Since $S\subseteq\mathcal X$, we have $\Psi_S\le\Psi$, with equality
and a common subgradient at every queried center inside the ball. The minimizer of $\Psi_S$ is
$\rho(v-x_S^*)$ with $x_S^*\in\mathrm{conv}(S)\subseteq\mathrm{conv}(\mathcal X)$, so it has the same
margin, and the subgradients of $\Psi_S$ obey the same bound on the ball. The same transcript is therefore
a valid run on $\Psi_S$, and some queried center $\bar p$ satisfies
$\min\Psi\le\Psi(\bar p)=\Psi_S(\bar p)\le\min\Psi_S+\varepsilon$.

The restricted primal
\[
P_S:=\min\Big\{\sum_{w\in S}\alpha_wf(w)+R_v\Big(\sum_{w\in S}\alpha_ww\Big):
\alpha\in\mathbb R_+^{S},\ \mathbf1^\top\alpha=1\Big\}
\]
equals $\frac{\rho}{2}\|v\|_2^2-\min\Psi_S$ by~\eqref{eq:prox_dual} applied to $f$ restricted to $S$, and
lies within $\varepsilon$ above the prox value. It is a rational convex quadratic program of polynomial
size, solvable exactly \citep{kozlov1979polynomial}. Its solution gives
$\hat x:=\sum_w\alpha_w^*w\in\mathrm{conv}(\mathcal X)$, and $f^{**}(\hat x)\le\sum_w\alpha_w^*f(w)$
yields the stated accuracy.

Restoring the agent index, use this approximation for the $x$-update~\eqref{eq:cadmm_x_step}, with
$v_i^k:=z_i^k-u_i^k$, at accuracy
$\varepsilon_k=\min\{1,\rho\tau_0^2/(2(k+1)^4)\}$ for any rational $\tau_0>0$. The prox errors then obey
$\|\hat x_i^k-\mathcal P_{i,\rho}(v_i^k)\|_2\le\tau_0/(k+1)^2$, which is summable, so the inexact iteration
converges \citep[Thm.~8]{eckstein1992douglas}, and every $x$-iterate, hence every ergodic average~\eqref{eq:ergodic_averages}, stays in
$\prod_i\mathrm{conv}(\mathcal X_i)$ as Corollary~\ref{cor:eps_stage1} requires. The rate~\eqref{eq:admm_ergodic_bounds} is proved for exact
proximal calls, and we do not claim it carries over unchanged to this inexact implementation.
\end{remark}

\begin{remark}[Separable shortcut]\label{rem:separable_prox}
On the canonical fixed-sum simplex domain, when
$f_i(x)=\sum_jg_{ij}(x^j)$ with discrete-convex coordinate costs, $f_i^{**}$ is the sum of the
coordinatewise chord interpolants. Its prox is then a separable convex resource-allocation problem with
one scalar multiplier, so the general ellipsoid construction is unnecessary. Separability alone does not
give this reduction on a general $M$-convex domain with additional base-polyhedron constraints.
\end{remark}

\section{Proofs of Proximity Bounds}\label{app:proximity_proofs}

We collect here the proof of Proposition~\ref{prop:proximity}. We retain the notation of Section~\ref{sec:proximity_bounds}: $G(x):=\sum_i f_i^{**}(x_i)$, $\tilde G(x):=\sum_i\tilde f_i(x_i)$, $\mathcal{F}:=\big\{x\in\prod_i\mathrm{conv}(\mathcal{X}_i):\sum_i A_i x_i=b\big\}$, and $\Delta:=\sum_i\Delta_i$ as in~\eqref{eq:envelope_gap_total}.

\paragraph{Proof of Proposition~\ref{prop:proximity}.}
We assemble four ingredients and combine them.

\emph{(i) Quadratic growth at $\tilde x^*$.} $\tilde G$ is $\mu$-strongly convex on $\mathcal{F}$ and $\tilde x^*\in\mathcal{F}$ is its constrained minimizer. Fix $t\in(0,1]$ and set $x_t:=(1-t)\tilde x^*+tx^*$, which lies in $\mathcal{F}$ because $\mathcal{F}$ is convex and $x^*,\tilde x^*\in\mathcal{F}$. Strong convexity along the segment gives $\tilde G(x_t)\le(1-t)\tilde G(\tilde x^*)+t\tilde G(x^*)-\tfrac{\mu}{2}t(1-t)\|x^*-\tilde x^*\|_2^2$, while minimality gives $\tilde G(\tilde x^*)\le\tilde G(x_t)$. Chaining the two, cancelling $(1-t)\tilde G(\tilde x^*)$ and dividing by $t>0$ leaves $\tilde G(\tilde x^*)\le\tilde G(x^*)-\tfrac{\mu}{2}(1-t)\|x^*-\tilde x^*\|_2^2$; letting $t\downarrow 0$ gives
\begin{align}
\tfrac{\mu}{2}\|x^*-\tilde x^*\|_2^2 \;\le\; \tilde G(x^*) - \tilde G(\tilde x^*).\label{eq:prox_step1}
\end{align}

\emph{(ii) Pointwise envelope inequality.} The section hypotheses give $\tilde G\le G$ on $\prod_i\mathrm{conv}(\mathcal{X}_i)$, while~\eqref{eq:envelope_gap_total} yields $G(x)-\tilde G(x)=\sum_i[f_i^{**}(x_i)-\tilde f_i(x_i)]\le\sum_i\Delta_i=\Delta$ on the same set. Hence on $\mathcal{F}$,
\begin{align}
\tilde G(x^*)\le G(x^*),\qquad \tilde G(\tilde x^*)\ge G(\tilde x^*) - \Delta.\label{eq:prox_step2}
\end{align}

\emph{(iii) LCE optimality.} Since $\tilde x^*\in\mathcal{F}$ and $x^*$ minimizes $G$ on $\mathcal{F}$,
\begin{align}
G(x^*) \le G(\tilde x^*).\label{eq:prox_step3}
\end{align}

\emph{(iv) Combine.} Substituting~\eqref{eq:prox_step2} and~\eqref{eq:prox_step3} into~\eqref{eq:prox_step1},
\begin{align*}
\tfrac{\mu}{2}\|x^*-\tilde x^*\|_2^2 \;\le\; \tilde G(x^*) - \tilde G(\tilde x^*) \;\le\; G(x^*) - \big(G(\tilde x^*)-\Delta\big) \;\le\; \Delta,
\end{align*}
and rearranging yields~\eqref{eq:proximity}.\qed

\section{Experimental Setup and Data Generation}\label{app:experimental_setup}

This appendix details the data-generation protocol, the implementation of the two-stage algorithm, the instance construction and solver settings, and the benchmark formulations used in Section~\ref{sec:numerical_illustration}.

\subsection{Synthetic Data Generation}\label{app:data_generation}

Each instance is parameterized by the number of regions $n$, the number of products $N$, and a random seed. Per-product data are drawn independently from distributions whose parameters were fit to a shadow dataset; the fitted values are the cross-tabulated medians of that dataset and are held fixed across all $(n,N)$.

\paragraph{Received quantities.} The quantity of product $i$ is lognormal on the positives, $\log q_i\sim\mathcal{N}(3.56,\,1.49^2)$, rounded to the nearest integer, with $q_i=0$ with probability $4\times10^{-4}$ (products with $q_i=0$ carry no mass and are discarded, since the cost~\eqref{eq:placement_cost} is then degenerate).

\paragraph{Demand and demand shares.} For each product the regional demand shares $(\rho_i^1,\ldots,\rho_i^n)$ are drawn from a symmetric Dirichlet distribution $\mathrm{Dir}(\alpha_0/n,\ldots,\alpha_0/n)$ with concentration $\alpha_0=60$; this simulates a per-region mean ($1/n$) and coefficient of variation observed in an empirical dataset. The demand target is $\mu_i^j=\rho_i^j\,(\,\sum_{j'}y_i^{j'}+q_i)$, the share applied to the total post-receipt inventory.

\paragraph{On-hand inventory.} On-hand inventory is sparse. With probability $0.16$ a product has no on-hand stock in any region. Otherwise the total on-hand $\sum_j y_i^j$ is lognormal, $\log(\sum_j y_i^j)\sim\mathcal{N}(2.70,\,1.66^2)$, distributed across regions by a Bernoulli$(0.36)$ activity mask (resampled if empty) followed by a uniform (Dirichlet$(1,\ldots,1)$) split over the active regions, rounded to integers.

\paragraph{Cost.} Given $(q_i,\mu_i^j,y_i^j)$, the placement cost~\eqref{eq:placement_cost} is assembled in the separable-quadratic form $f_i(x)=\sum_j\tfrac12 d_i^j (x^j)^2+c_i^j x^j+\mathrm{const}_i$ with $d_i^j=2/\mu_i^j$, $c_i^j=-2(\mu_i^j-y_i^j)/\mu_i^j$, and an additive constant collecting the $y_i^j$ terms. Since $\mu_i^j>0$, every $d_i^j>0$, so each $f_i$ is strictly convex and serves directly as the strongly convex parent surrogate of Section~\ref{sec:strongly_convex_extension}.

\paragraph{Weights and capacities.} Each product is assigned an integer arc-resource weight $w_i$ drawn uniformly from $\{1,\ldots,10\}$, giving the per-agent diagonal coupling block $A_i=w_i I$. The arc capacities are an even split of the total weighted mass, $b^j=\lfloor(\sum_i w_i q_i)/n\rfloor$ for $j<n$ with the remainder assigned to $j=n$, so that $\sum_j b^j=\sum_i w_i q_i$, which is necessary, though not sufficient, for an integer point to satisfy the coupling constraint $\sum_i w_i x_i^j=b^j$. With weights $w_i>1$ the stacked matrix $[A_1\cdots A_N]$ is not totally unimodular, so integrality of the relaxation is not guaranteed; it was fractional on every instance in the sweep. Integer feasible points nonetheless exist throughout: the exact MILP found one on every cell it was run on.

\subsection{Implementation of the Two-Stage Algorithm}\label{app:admm_impl}

This subsection specializes the generic Stage~1 solver of Appendix~\ref{app:admm_details} to the separable cost~\eqref{eq:placement_cost} and gives the coupling projection of Section~\ref{sec:num_solvers}.

\paragraph{Proximal step on the envelopes.}
For the separable cost of Section~\ref{sec:numerical_illustration}, write $g_i^j(t)=\tfrac12d_i^jt^2+c_i^jt$ for coordinate $j$ of agent $i$ (constants dropped) and $\hat g_i^j$ for its chord interpolant at the integers, so that $f_i^{**}(x)=\sum_j\hat g_i^j(x^j)$ on $\mathrm{conv}(\mathcal X_i)$ by Proposition~\ref{rem:separable_decomposition}. On the cell $[k,k+1]$, $\hat g_i^j$ is affine with slope $s_k=g_i^j(k+1)-g_i^j(k)=d_i^j(k+\tfrac12)+c_i^j$, increasing in $k$. The proximal step~\eqref{eq:cadmm_x_step} at $v=z_i^k-u_i^k$ is
\[
\min_{x\ge0}\ \sum_j\Big[\hat g_i^j(x^j)+\tfrac\rho2(x^j-v^j)^2\Big]\quad\text{s.t.}\quad\textstyle\sum_jx^j=q_i .
\]
For a multiplier $\nu$ of the sum constraint the problem separates. Coordinate $j$ minimizes the convex piecewise quadratic $\hat g_i^j(t)+\nu t+\tfrac\rho2(t-v^j)^2$ over $t\ge0$, whose minimizer on cell $k$ is $t_k=v^j-(s_k+\nu)/\rho$. The condition $t_k\ge k$ reads $a\ge k(1+d_i^j/\rho)$ with $a:=v^j-(\tfrac12d_i^j+c_i^j+\nu)/\rho$, so the optimal cell is $k^*=\lfloor a/(1+d_i^j/\rho)\rfloor$ and
\[
x^j(\nu)=\max\Big\{0,\ \min\big\{k^*+1,\ \max\{k^*,\,t_{k^*}\}\big\}\Big\}.
\]
The map $\nu\mapsto\sum_jx^j(\nu)$ is continuous, nonincreasing and piecewise linear, so the multiplier with $\sum_jx^j(\nu)=q_i$ is found by a bracketed one-dimensional search. The bracket is centred at the multiplier of the same problem with $g_i^j$ in place of $\hat g_i^j$, which is a closed-form weighted-simplex projection and lies within $O(\max_jd_i^j)$ of the root, and the search alternates bisection with an exact Newton step on the coordinates interior to their cell. Each evaluation costs $O(n)$ per agent and is vectorized over all agents. The $z$-step uses the diagonal closed form of Appendix~\ref{app:admm_details} with $a_i=w_i\mathbf 1$.

\paragraph{Projection onto the coupling.}
Given the Stage~1 iterate $x^*$ and per-agent boxes $\ell_i\le u_i$, the coupling projection of Section~\ref{sec:num_solvers} solves
\[
\min_x\ \tfrac12\sum_i\|x_i-x_i^*\|_2^2\quad\text{s.t.}\quad\sum_iA_ix_i=b,\qquad x_i\in S_i:=\{x:\mathbf 1^\top x=q_i,\ \ell_i\le x\le u_i\}.
\]
Dualizing the coupling with $\eta\in\mathbb R^m$, each agent's best response is the projection $x_i(\eta)=\Pi_{S_i}(x_i^*-A_i^\top\eta)$ onto a boxed simplex, computed by the same kind of vectorized one-dimensional search. The residual $g(\eta):=\sum_iA_ix_i(\eta)-b$ is the gradient of the concave dual function, hence monotone (nonincreasing) and piecewise affine, with generalized Jacobian $-H(\eta)$, $H=\sum_iA_iP_iA_i^\top$, where $P_i$ is the projector $I_{F_i}-\mathbf 1_{F_i}\mathbf 1_{F_i}^\top/|F_i|$ onto the coordinates $F_i$ of $x_i(\eta)$ strictly inside their box. We apply a damped semismooth Newton method from $\eta=0$: the direction is $H^{\dagger}g(\eta)$, the step is chosen by Armijo backtracking on the dual function, and the iteration stops when $\|g\|_\infty\le10^{-9}\max(1,\|b\|_\infty)$. Each iteration costs one projection per agent, one pass over the agents to assemble $H$ ($O(Nn^2)$ for the diagonal blocks $A_i=w_iI$ used here) and an $m\times m$ solve, so the projection is negligible next to one ADMM iteration whenever it takes few steps, which is the case here: on the comparison grid it took at most four Newton steps and $0.3$\,s, and on the scaling runs at most $30$\,s except for one run at $118$\,s. After projection the coupling residual $\|\sum_iA_i\tilde x_i-b\|_\infty$ was below $0.04$ on every run. With $\ell_i=\lfloor x_i^*\rfloor$ and $u_i=\lceil x_i^*\rceil$ (integral coordinates are fixed) the problem may be infeasible when the cells cannot absorb the miss; the implementation then repeats the projection with $\ell_i=0$, $u_i=q_i\mathbf 1$, which is feasible whenever~\eqref{eq:convex_relaxation} is.

\subsection{Instance Construction and Solver Settings}\label{app:solver_settings}

For each $(n,N,\text{seed})$ we draw functions until $N$ valid products remain (after discarding degenerate $q_i=0$ draws) and take the first $N$; the construction is deterministic in the triple, so every method sees the identical instance. The comparison grid is $n\in\{4,8,16,32\}$, $N\in\{2^8,\ldots,2^{14}\}$ and seeds $\{0,\ldots,9\}$. The scaling runs of Table~\ref{tab:scale} use the same seeds at $N\in\{2^{13},\ldots,2^{18}\}$ and seeds $\{0,1,2\}$ at $N\in\{2^{19},2^{20}\}$.

\paragraph{Hardware.} All experiments ran on one machine with an Intel Xeon Platinum 8259CL at $2.5$\,GHz (16 cores, 32 hardware threads) and $124$\,GiB of visible memory. All times are wall-clock seconds on this machine.

\paragraph{Two-stage settings.} Stage~1 is consensus ADMM (Appendix~\ref{app:admm_details}) on the envelopes $f_i^{**}$ with initial penalty $\rho=1$, over-relaxation $1.5$, residual balancing of $\rho$ (rescaling when the primal and dual residuals differ by more than a factor $5$, factor capped at $100$, frozen after $500$ iterations), absolute tolerance $10^{-7}$, relative tolerance $10^{-5}$ and at most $1000$ iterations. The relative tests scale with the norm of the iterates, which grows with $N$; the coupling projection of Section~\ref{sec:num_solvers} and Appendix~\ref{app:admm_impl} removes the resulting coupling miss, so the tolerance is the same at every $N$. Stage~2(a) enumerates $\mathcal X_i\cap\mathcal N(\tilde x_i)$ when $\binom{|I_i|}{p_i}\le24$ and computes the weights $\alpha_i$ by a per-agent LP, and otherwise runs column generation on~\eqref{eq:decomposition_lp} from the offset-sweep basis of Appendix~\ref{app:offset_sweep}, over a pool of worker processes (for this separable cost the sweep basis alone is already a value-preserving decomposition; see Section~\ref{sec:num_solvers}). Before Stage~2(a), any residual round-off off the slice $\sum_jx_i^j=q_i$ is removed by subtracting the excess in proportion to the coordinates. Stage~2(b) is Algorithm~\ref{alg:purification} with $\hat m=m-1$, the block rule of Section~\ref{sec:num_solvers}, a support threshold of $10^{-9}$ on the weights, and a check that the starting weights satisfy the coupling to within $10^{-6}\max(1,\|b\|_\infty)$.

\paragraph{Comparison protocol.} The comparison grid ran as four concurrent jobs, one per $n$, each pinned with \texttt{taskset} to four physical cores and their hyperthreads ($8$ hardware threads), the same set for every method run for that $n$; SCIP, whose branch and bound is sequential, ran on one core per $n$.

\subsection{Benchmarks}\label{app:benchmarks}

The natural exact benchmark is the MIQP that states~\eqref{eq:base_problem} directly with integer variables and the quadratic cost. Both MIQP solvers we tried struggled on these instances: Xpress did not close the gap even on the smallest cells, and SCIP was the slowest method on every cell it certified. We therefore also solve two piecewise-linear MILP reformulations, which fared markedly better.

All benchmarks start cold, receive no information from the two-stage method, and are timed from the start of model construction to the end of the solve. A cell is certified when the solver reports a relative gap of at most $10^{-4}$ between incumbent and bound on the true cost, so ``certified'' means optimal to within that tolerance; the cost constant is passed to each solver as an objective offset so that its relative gap is measured on the true cost. Each run has a $3600$\,s limit. If any seed of a cell $(n,N)$ reaches the limit without certifying, the method is not run on cells $(n',N')\ne(n,N)$ with $n'\ge n$ and $N'\ge N$; the remaining seeds of $(n,N)$ still run. Every MILP run records its incumbent and bound at each new incumbent, which gives the time to the first feasible solution and to the first solution at least as good as the two-stage one.

\paragraph{Exact piecewise-linear MILP.} Each coordinate cost $g_i^j(t)=\tfrac12 d_i^j t^2+c_i^j t$ is replaced by its integer-point piecewise-linear interpolant. Split $x_i^j\in[0,q_i]$ into unit increments $\omega_{ijk}\in[0,1]$, $k=0,\ldots,q_i-1$, with marginal cost $s_{ijk}=g_i^j(k+1)-g_i^j(k)=d_i^j(k+\tfrac12)+c_i^j$, so $x_i^j=\sum_k \omega_{ijk}$ and the cost is the linear $\sum_{ijk}s_{ijk}\omega_{ijk}$ plus a constant. By convexity the slopes increase in $k$, so for an integer value of the aggregate $x_i^j$ an optimal choice of increments takes values in $\{0,1\}$ and fills them cheapest-first, and the linear cost equals $g_i^j(x_i^j)$. It therefore suffices to impose integrality on the $Nn$ aggregates rather than on the $n\sum_iq_i$ increments, with no SOS2 or ordering constraints, and the model is an exact integer reformulation of~\eqref{eq:base_problem}. The constraints are $\sum_jx_i^j=q_i$ and $\sum_iw_ix_i^j=b^j$. We solve it with FICO Xpress $9.8$ through its Python API with $8$ threads, \texttt{miprelstop} $10^{-4}$, \texttt{mipabsstop} $10^{-9}$, the default cutoff slacks \texttt{mipaddcutoff} and \texttt{miprelcutoff} set to $0$ (by default they discard nodes within about $10^{-5}$ of the incumbent), and \texttt{maxtime} $3600$. With the identity coupling $A_i=I$ the constraint matrix would be that of a bipartite transportation problem, totally unimodular, and the LP relaxation would already be integral; the weights $w_i>1$ break this, the LP relaxation was fractional on every instance, and the solver branches.

\paragraph{Relaxed-coupling MILP.} The same model and settings, with each coupling row replaced by the range $b^j-r\le\sum_iw_ix_i^j\le b^j+r$, where $r=\|\sum_iw_i\bar x_i-b\|_\infty$ is the residual of the two-stage solution on the same instance.

\paragraph{SCIP MIQP.} Integer variables $x_i^j\in\{0,\ldots,q_i\}$, the two families of linear equalities, and one epigraph variable $t_i^j\ge g_i^j(x_i^j)$ per term, with objective $\sum_{ij}t_i^j$ plus the constant as an offset. We solve it with SCIP $10.0$ through PySCIPOpt with \texttt{limits/gap} $10^{-4}$, \texttt{limits/absgap} $10^{-9}$, \texttt{limits/time} $3600$ and default settings otherwise. SCIP's branch and bound runs in one thread. With default settings, SCIP~$10.0$ aborted with a heap-corruption error at the first branching node on most instances with $n\ge8$; enabling its primal heuristics one at a time traced the abort to the NLP diving heuristic, and with that heuristic disabled (\texttt{heuristics/nlpdiving/freq} $=-1$) the same instances solve normally, so it is disabled in all reported runs. As a safeguard each SCIP solve runs in its own process; a seed on which SCIP still aborted would be recorded as a solver failure, not a timeout, and excluded from the times. Reported SCIP times are model construction plus solve, without the process start-up.

\paragraph{The plain MIQP in Xpress.} Xpress on the quadratic-objective MIQP finds good incumbents quickly but does not close the gap: its continuous relaxation is the QP relaxation, which lies below the lower convex envelope by up to $d_i^j/8=1/(4\mu_i^j)$ for every fractional coordinate, and branching closes this slowly on these degenerate instances. In our tests it stopped at relative gaps between $7\times10^{-4}$ and $1.5\times10^{-3}$ after $300$ to $900$\,s on the cells with $n=4$ and $N\le2^{10}$, which the other methods certify in seconds, and relaxing the coupling by $r$ left the gap essentially unchanged. We therefore do not report it.

\section{Hypersimplex Offset-Sweep Decomposition}\label{app:offset_sweep}

When $\binom{|I_i|}{p_i}$ is too large to enumerate all hypersimplex vertices above $\lfloor x_i^*\rfloor$, the practical column-generation approach of Remark~\ref{rem:practical_vs_theoretical} becomes the method of choice: a revised simplex loop that maintains a basis of at most $n$ vertices, the rank of the decomposition LP once the redundant normalization row is removed, and uses $M$-convex pricing to introduce new vertices as needed. This approach requires an \emph{initial feasible basis}, a set of at most $n$ hypersimplex vertices whose convex hull contains $x_i^*$ (equivalently, whose weights $\{\alpha_l\}$ can be chosen nonnegatively with $\sum_l\alpha_l=1$ and $\sum_l\alpha_l(\lfloor x_i^*\rfloor+e_{S_l})=x_i^*$). For a separable cost on the canonical domain the sweep below is itself the full Stage~2(a) (Proposition~\ref{rem:separable_decomposition}).

Producing such an initial decomposition is a classical problem: it asks to decompose a point in the hypersimplex $\Delta(p_i,n)=\{\theta\in[0,1]^n:\sum_j\theta_j=p_i\}$ into a convex combination of $0/1$ vertices. The algorithm is closely related to the \emph{dependent rounding} schemes of \citet{gandhi2006dependent,chekuri2010dependent}.

We describe one simple deterministic implementation, an \emph{offset-sweep} procedure running in $\mathcal{O}(n^2)$ time and producing $K\le n$ vertices. The resulting decomposition is purely geometric (it reproduces $x_i^*$ via its marginals) and in general does \emph{not} achieve the value identity~\eqref{eq:value_decomposition}. It nonetheless serves as a valid warm-start basis for column generation, which then refines the cost via $M$-convex pricing. The vertices it returns lie in $\lfloor x_i^*\rfloor+\{e_S:|S|=p_i\}$, hence in $\mathcal{X}_i$ for the canonical domain of Section~\ref{sec:mconvexity}. For a general $M$-convex domain they need not be feasible and must be filtered against the value oracle, and the two things the sweep is used for then separate. A single surviving pattern already supplies the initial point $v^0$ required by Theorem~\ref{thm:poly_extension}, and if none survives the phase-one minimization~\eqref{eq:local_phase_one} supplies it; a feasible column-generation basis is more demanding, since it needs a surviving subfamily whose convex hull still contains $x_i^*$, and filtering can destroy that even when patterns survive. In that case one starts column generation from the basis produced by the ellipsoid route of Theorem~\ref{thm:poly_extension} and Remark~\ref{rem:poly_decomposition}, which needs only $v^0$.

\paragraph{Setup.} Fix an agent $i$ and drop the subscript. Let $\theta:=x^*-\lfloor x^*\rfloor\in[0,1)^n$ be the fractional part, with $\sum_j\theta_j=p\in\mathbb{N}$. The goal is to write $\theta=\sum_{l=1}^K\alpha_l e_{S_l}$ as a convex combination of hypersimplex vertices $e_{S_l}$ with $|S_l|=p$ and $K\le n$, so that $x^*=\sum_l\alpha_l(\lfloor x^*\rfloor+e_{S_l})$.

\paragraph{Idea.} Choose an ordering $\pi$ of the coordinates and compute the cumulative sums $C_k=\sum_{l=1}^k\theta_{\pi(l)}$. For any offset $t\in[0,1)$, define the rounding pattern $S(t)\subseteq\{1,\ldots,n\}$ by $\pi(k)\in S(t)$ iff $\lfloor C_{k-1}+t\rfloor<\lfloor C_k+t\rfloor$ (i.e., coordinate $\pi(k)$ is rounded up when a floor is crossed). As $t$ sweeps through $[0,1)$, the pattern $S(t)$ changes exactly when $t$ crosses one of the breakpoints $t_k=(1-\{C_k\})\bmod 1$, producing at most $n$ distinct patterns. Weighting each pattern by the length of its interval yields the desired decomposition.

\begin{algorithm}[h!]
\SetAlgoNoLine \small
\textbf{Input:} Fractional part $\theta \in [0,1)^n$ with $\sum_j \theta_j = p \in \mathbb{N}$\\[4pt]
Choose an ordering $\pi$ of $\{1,\ldots,n\}$\;
Compute cumulative sums: $C_0 \leftarrow 0$;\; $C_k \leftarrow C_{k-1} + \theta_{\pi(k)}$ for $k=1,\ldots,n$\;
Compute breakpoints: $t_k \leftarrow (1 - \{C_k\}) \bmod 1$ for $k=0,\ldots,n-1$, where $\{C_k\} = C_k - \lfloor C_k \rfloor$\;
Sort the distinct breakpoints: $0 = t_{(1)} < t_{(2)} < \cdots < t_{(K)}$, $t_{(K+1)} = 1$\;
\For{$l = 1,\ldots,K$}{
    Weight: $\alpha_l \leftarrow t_{(l+1)} - t_{(l)}$\;
    Pick $\bar{t} \in (t_{(l)}, t_{(l+1)})$; vertex $e_{S_l}$: $\pi(k) \in S_l$ iff $\lfloor C_{k-1} + \bar{t}\rfloor < \lfloor C_k + \bar{t}\rfloor$\;
}
\textbf{Output:} Vertices $\{e_{S_l}\}$ and weights $\{\alpha_l\}$ with $\theta = \sum_{l=1}^K \alpha_l e_{S_l}$, $K \leq n$
\caption{Hypersimplex Decomposition via Offset Sweep}
\label{alg:hypersimplex}
\end{algorithm}

\paragraph{Correctness and complexity.} The algorithm produces a valid geometric decomposition: $\theta=\sum_l\alpha_l e_{S_l}$ with $\sum_l\alpha_l=1$, $\alpha_l\ge 0$, and $|S_l|=p$ for every $l$ (the last property follows from the hyperplane constraint $\sum_j\theta_j=p$, which forces exactly $p$ floor crossings for every value of $t$). The marginal condition $\sum_{l:\,j\in S_l}\alpha_l=\theta_j$ holds for each coordinate $j$. Runtime is $\mathcal{O}(n^2)$: $\mathcal{O}(n)$ for cumulative sums and breakpoints, and $\mathcal{O}(n)$ per interval to compute $S_l$.

\paragraph{Cost-obliviousness.} The decomposition depends only on the geometry of $\theta$, not on the cost function $f$. Consequently the value identity~\eqref{eq:value_decomposition} does \emph{not} hold in general: $\sum_l\alpha_l f(\lfloor x^*\rfloor+e_{S_l})$ may exceed $f^{**}(x^*)$. The exception is a separable cost
on the canonical domain, for which every local decomposition, and hence the sweep, is value-preserving
(Remark~\ref{rem:separable_decomposition}). In general this is expected and acceptable: the algorithm is used as a warm-start basis for the column-generation loop of Remark~\ref{rem:practical_vs_theoretical}, which then refines the cost via $M$-convex pricing. Running the offset sweep with multiple orderings $\pi$ produces richer candidate bases, which can reduce the number of pricing iterations.

\end{document}